\documentclass[12pt,draft]{article}
\usepackage{amsthm,amssymb,amsmath,amscd,amstext,mathrsfs}
\usepackage{latexsym}

\newtheorem{theorem}{Theorem}[section]
\newtheorem{definition}{Definition}[section]
\newtheorem{lemma}{Lemma}[section]

\newtheorem*{remark}{{\it Remark}}

\newcommand{\nc}{\newcommand}
\nc{\C}{{\mathbb C}}
\nc{\R}{{\mathbb R}}
\nc{\HH}{{\mathbb H}}
\nc{\Z}{{\mathbb Z}}
\nc{\N}{{\mathbb N}}
\nc{\dd}{{\rm d}}
\nc{\DD}{{\bf d}}
\nc{\ii}{{\bf i}}
\nc{\jj}{{\bf j}}
\nc{\kk}{{\bf k}}
\nc{\co}{{\mathscr O}}
\nc{\ck}{{\mathscr K}}
\nc{\cj}{{\mathscr J}}
\nc{\cm}{{\mathscr M}}
\nc{\crr}{{\mathscr R}}
\nc{\cs}{{\mathscr S}}
\nc{\ce}{{\mathscr E}}
\nc{\cf}{{\mathscr F}}
\nc{\cg}{{\mathscr G}}
\nc{\ZZ}{\overline{Z}}
\nc{\tr}{\mathop{\rm tr}\nolimits}
\nc{\ext}{\mathop{\rm Ext}\nolimits}
\nc{\su}{{\mathfrak s}{\mathfrak u}(2)}
\nc{\so}{{\mathfrak s}{\mathfrak o} (4)}

\begin{document} 

\title{Harmonic functions and instanton moduli spaces on the 
multi-Taub--NUT space} 

\author{G\'abor Etesi\\
\small{{\it Department of Geometry, Mathematical Institute, Faculty of 
Science,}}\\
\small{{\it Budapest University of Technology and Economics,}}\\
\small{{\it Egry J. u. 1, H \'ep., H-1111 Budapest, Hungary
\footnote{e-mail: {\tt etesi@math.bme.hu}}}}\\
Szil\'ard Szab\'o\\
\small{{\it Department of Geometry, Mathematical Institute, Faculty of
Science,}}\\
\small{{\it Budapest University of Technology and Economics,}}\\
\small{{\it Egry J. u. 1, H \'ep., H-1111 Budapest, Hungary
\footnote{e-mail: {\tt szabosz@math.bme.hu}}}}}

\maketitle

\pagestyle{myheadings}
\markright{G. Etesi, Sz. Szab\'o: Instantons over the multi-Taub--NUT space}  

\thispagestyle{empty}

\begin{abstract}
Explicit construction of the basic SU$(2)$
anti-instantons over the multi-Taub--NUT geometry via the
classical conformal rescaling method is exhibited. These anti-instantons
satisfy the so-called weak holonomy condition at infinity with respect 
to the trivial flat connection and decay rapidly. The resulting unit 
energy anti-instantons have trivial holonomy at infinity.

We also fully describe their unframed moduli space and find that it 
is a five dimensional space admitting a singular disk-fibration over 
$\R^3$.

On the way, we work out in detail the twistor space of the multi-Taub--NUT
geometry together with its real structure and 
transform our anti-instantons into holomorphic vector bundles over the 
twistor space. In this picture we are able to demonstrate that our
construction is complete in the sense that we have constructed a full 
connected component of the moduli space of solutions of the above type.

We also prove that anti-instantons with arbitrary high integer energy 
exist on the multi-Taub--NUT space.
\end{abstract}
\vspace{0.5cm}

\centerline{AMS Classification: Primary: 53C07; Secondary: 14F05, 53C28, 58J99}
\vspace{0.5cm}

\centerline{Keywords: {\it Multi-Taub--NUT space; Yang--Mills
instantons; Harmonic functions; Twistors}}


\section{Introduction}
\label{one}

The aim of this paper is to construct the most relevant anti-instanton 
moduli space over the multi-Taub--NUT geometry by elementary means.
 
An important class of non-compact but complete four
dimensional geometries is the collection of the so-called {\it 
asymptotically locally flat} (ALF) spaces including several mathematically 
as well as physically important examples. The flat $\R^3\times S^1$ 
plays a role in finite temperature Yang--Mills theories, 
the Euclidean Schwarzschild space \cite{haw} deals with quantum 
gravity and Hawking radiation. If the metric is additionally hyper-K\"ahler 
then the space is also called an {\it ALF gravitational instanton} 
(in the narrow sense). The flat $\R^3\times S^1$ is a straightforward 
example; non-trivial ones for this restricted class are provided by the 
multi-Taub--NUT (or $A_k$ ALF or ALF 
Gibbons--Hawking) spaces \cite{gib-haw} which carry supersymmetric 
solutions of string theory and supergravity models, the Atiyah--Hitchin 
manifold (and its universal double cover) \cite{ati-hit} describing the 
2-monopole moduli space over $\R^3$ and last but not least the recently 
constructed $D_k$ ALF spaces \cite{che-hit}.

The ALF asymptotics is a natural generalization of the well-known ALE 
(asymptotically locally Euclidean) one including the multi-Eguchi--Hanson 
geometries \cite{gib-haw}. Instanton theory over these
later spaces possessing several 
phenomena related with non-compactness (e.g. existence of four 
dimensional moduli spaces isometric to the original space)
is well-known due to the important paper of Kronheimer and Nakajima 
\cite{kro-nak} in which a full ADHM construction was established. The 
existence of this construction is in some sense not surprising 
because the original ADHM construction was designed for the flat 
$\R^4$ and all ALE spaces arise by an algebro-geometric deformation of 
the {\it flat} quotients $\C^2/\Gamma$ where $\Gamma\subset {\rm SU}(2)$ 
are various discrete subgroups \cite{kro}.

The natural question arises: what about instanton 
theory over ALF spaces? Unlike the ALE geometries, these are 
essentially non-flat spaces in the sense above therefore one may expect that 
instanton theory somewhat deviates from that over the flat $\R^4$. Some 
general questions have been answered recently 
\cite{ete-jar} and these investigations pointed out that in spite of 
their more transcendental nature, open Riemannian spaces with ALF asymptotics 
rather resemble compact four manifolds at least from the point of view of 
instanton theory. To test this interesting observation more carefully in 
this paper we work out the simplest moduli space over the 
multi-Taub--NUT geometry: we will find that this moduli space is five 
dimensional furthermore its part containing concentrated instantons 
looks like a ``collar'' of the original manifold, supporting the analogy 
with the compact case. We note that this apparent compactness is related to the 
existence of a smooth compactification of ALF spaces motivated by $L^2$ 
cohomology theory \cite{hau-hun-maz}.

The paper is organized as follows. In Sect. \ref{two} we quickly 
summarize two important tools for solving the self-duality equations 
over an anti-half-conformally flat space: the conformal rescaling method 
of Jackiw--Nohl--Rebbi \cite{jac-noh-reb,ati-hit-sin} as well as the 
Atiyah--Ward correspondence \cite{ati-war}. Via conformal rescaling one 
constructs instantons out of positive harmonic functions with at most 
pointlike singularities and appropriate decay toward infinity while the 
Atiyah--Ward correspondence establishes a link between instantons and 
holomorphic vector bundles. Since these methods meet in 
Penrose' twistor theory  we shall also briefly outline it here.

Then we introduce ALF spaces. Referring to recent results on 
instanton theory over these spaces \cite{ete-jar} we precisely define the 
class of anti-instantons which are expected to form nice moduli spaces. 
These solutions must obey two conditions: the so-called 
{\it weak holonomy condition} with respect to some smooth flat connection and 
the {\it rapid decay condition} (cf. Definition \ref{feltetelek} here). The 
former condition, albeit looks like an analytical one, in fact deals 
with the topology of infinity of the ALF space only 
\cite[Theorem 2.3]{ete-jar} meanwhile the latter 
one controls the fall-off properties of an anti-instanton hence is indeed 
analytical in its nature.

In Sect. \ref{three}, taking probably the most relevant ALF example namely 
the multi-Taub--NUT series, following 
Andreson--Kronheimer--LeBrun \cite{and-kro-leb} and LeBrun \cite{leb} 
we review this geometry with special attention to the description of its 
generic complex structures as algebraic 
surfaces in $\C^3$. We also identify a finite number of distinguished 
complex structures that can be realized as blowups of these algebraic 
surfaces. The group $S^1$ acts on these spaces via isometries and this 
action has isolated fixed points, called NUTs. 

Then in Sect. \ref{four} we come to our main results. We shall prove via 
conformal rescaling that over the multi-Taub--NUT space 
$(M_V, g_V)$ unframed $L^2$ moduli spaces of 
SU(2) anti-instantons obeying both the aforementioned weak holonomy 
condition with respect to the trivial flat connection $\nabla_\Theta$ as 
well as the rapid decay condition are not empty. Since the input in this 
construction is a positive harmonic function, the so-called {\it 
non-parabolic manifolds} \cite{gri, heb, var} are the key concept here 
which admit an abundance of {\it positive
minimal Green functions} $G(\cdot ,y)$ with arbitrarily prescribed 
pointlike singularity $y\in M_V$ (cf. Definition \ref{green} here).
Non-parabolicity deals with the volume growth toward infinity and is
closely related to the ALF property of our space as it was 
observed by Minerbe \cite[Theorem 0.1]{min1}. These functions provide us 
with a five real paremeter family
\[f_{y,\lambda}=1+\lambda G(\cdot ,y)\]
of positive harmonic functions hence that of non-gauge equivalent 
irreducible SU(2)-anti-self-dual connections of unit 
energy on the positive chiral spinor bundle---parameterized by a
point $y\in M_V$ and a ``concentration parameter'' $\lambda\in
(0, +\infty )$. These solutions obey the weak holonomy condition with respect 
to the trivial flat connection as well as decay rapidly. This way we 
obtain a familiar ``collar'' for the original manifold 
(cf. Theorem \ref{galler} here). Since we know from 
\cite[Theorem 3.2]{ete-jar} that the unit energy unframed moduli space of 
these irreducible solutions is a 
smooth five dimensional manifold it follows that this collar in fact 
realizes an open subset of it (or one of its connected component) and 
describes the regime of concentrated 
unit energy anti-instantons with trivial holonomy at infinity.

Then we proceed further and construct additional anti-instantons in the 
``centerless'' direction i.e., when the concentration parameter 
$\lambda$ becomes infinite. These solutions just correspond to pure 
Green functions. However in this limit anti-instantons corresponding to the 
NUTs are reducible to U$(1)$ consequently it is natural to talk about (a 
connected component of) the unframed 
moduli space $\widehat{\cm}(1,\Theta )$ including these reducible 
solutions as well. Therefore we will find in particular that this 
extended moduli space is not a manifold anymore rather admits usual 
singularities in the reducible points moreover that our moduli space 
itself admits a singular fibration 
\[\Phi :\widehat{\cm}(1,\Theta )\longrightarrow\R^3\]
whose generic fibers are $2$ dimensional open while the singular ones 
are $1$ dimensional semi-open balls (cf. Theorem \ref{moduluster} here). 
The aforementioned anti-instantons with infinite concentration parameter 
$\lambda =+\infty$ are represented in the generic case by the centers 
of these open $2$-balls meanwhile in the singular case by the closed end 
of these semi-open $1$-balls. For a given generic fiber 
the $\lambda ={\rm const.}<+\infty$ circles around the center 
describe anti-instantons ``centered'' about points in $(M_V,g_V)$ along 
the same isometry orbit and with fixed finite concentration parameter. 
The singular fibers correspond to the fixed points of the 
isometric $S^1$-action hence in a given singular fiber there is a 
unique anti-instanton with $\lambda ={\rm const.}\leq +\infty$. In 
other words, all ``centerless'' anti-instantons stemming from pure Green 
functions with pointlike singularities along a given 
isometry orbit are gauge equivalent.

We will also be able to demonstrate that similar moduli spaces
with arbitrary high energy and trivial holonomy at infinity are 
non-empty according to previous expectations.

For a comparison we also construct the moduli space of 
rapidly decaying unit energy SU(2) instantons with trivial holonomy in 
infinity over the flat $\R^3\times S^1$ (cf. Theorem \ref{kaloron} here). 

In the last section, Sect. \ref{five} we round things off by proving, as 
already mentioned, that our construction is complete in the sense that 
at least a whole connected component of the basic moduli space emerges 
this way. This will require however a quite tedious work.

To this end first we provide a detailed construction of the twistor space of 
the multi-Taub--NUT geometry together with its 
real structure {\it \`a la} Hitchin \cite{bes,hit1,hit3}. 
Then we transform our anti-instantons constructed out of harmonic 
functions into holomorphic vector bundles on the twistor space in the 
spirit of the Atiyah--Ward correspondence \cite{ati-war}. More precisely, 
by the aid of Atiyah \cite{ati} we work out a smooth 
compactification of the twistor space which makes it possible to 
identify the harmonic functions used to construct anti-instantons 
with elements of certain Ext groups of the 
compactified twistor space (cf. Theorem \ref{harmonikus-ext} here). 
A novelty of our compactification compared to Atiyah's 
is that our resulting Ext groups have the correct dimension (cf. Lemma 
\ref{ext-lemma}) i.e., they inculde the constant harmonic functions as 
well. (Note that without compactification the Ext groups would be infinite 
dimensional.)

These Ext groups can be used to obtain holomorphic bundles on the twistor 
space. In particular gluing data of the bundles 
corresponding to minimal Green functions can be worked out explicitly 
\cite{ati} and one can conclude that the isometric $S^1$-action induces an 
isomorphism of these bundles therefore the associated 
``centerless'' anti-instantons indeed must be gauge equivalent. This 
provides us with the picture of the moduli space mentioned above.

Working with the twistor space also can be viewed as a 
step toward establishing a general ADHM-like instanton factory 
\cite{che} which can be used to replace the conformal rescaling method; 
this is certainly necessary since our technique is not capable to grasp 
all higher energy moduli spaces. 

We conjecture that in fact we have obtained the full moduli 
space of unit energy, rapidly decaying anti-instantons with trivial 
holonomy at infinity in this paper i.e., the basic moduli space is 
connected.

Finally for a comparison we make a comment here on other 
gravitational instantons of ALG and 
ALH type introduced by Cherkis and Kapustin. By their definition, these 
geometries are parabolic hence one cannot expect the 
existence of positive minimal Green functions over them and accordingly, a 
canonical five parameter family of anti-instantons. Indeed, the simplest 
ALG space is the flat $\R^2\times T^2$ and it is known that the framed unit 
energy anti-instanton moduli space here is one dimensional only \cite{biq-jar}.
Hence in the ALE-ALH hierarchy apparently ALF spaces represent the 
closest analogues of compact four-manifolds from an instanton theoretic 
viewpoint. Analogously, in this hierarchy ALF spaces are the only 
ones whose $L^2$ cohomological compactification \cite{hau-hun-maz} is 
smooth.

\paragraph{Acknowledgement.} The second author would like to thank the 
Alfr\'ed R\'enyi Institute of Mathematics, where he had a scholarship 
during the early stages of this article, and Vincent Minerbe for useful 
discussions. Both authors were partially supported by OTKA grant No. 
NK81203 (Hungary).


\section{Constructing instantons}
\label{two}

In this section we quickly summarize two methods to obtain instantons 
over an anti-self-dual Riemannian manifold: the conformal rescaling 
method and the Atiyah--Ward correspondence. 
These two approaches meet each other in twistor theory therefore we also 
review the basic facts about it here.

Let us begin with the conformal rescaling method which 
was first used by Jackiw, Nohl and Rebbi to obtain instantons 
over the flat $\R^4$ \cite{jac-noh-reb}. Remember that an SU$(2)$ 
(anti-)self-dual connection $\nabla_A$ over a four dimensional oriented 
Riemannian manifold $(M,g)$ is a smooth SU$(2)$ connection on any 
SU$(2)$ rank two complex vector bundle $E$ over $M$ whose curvature satisfies 
the (anti-)self-duality equations $*_g F_A=\pm F_A$. The gauge 
equivalence class $[\nabla_A ]$ of this connection is 
called an {\it (anti-)instanton} if additionally its energy is 
finite: $\Vert F_A\Vert_{L^2(M)}<+\infty$. Note that this assumption 
on the energy is non-trivial only if $(M,g)$ is non-compact but complete. 
 
Let $(M,g)$ be an oriented Riemannian spin four-manifold, possibly 
non-compact but complete. Take the Levi--Civita connection $\nabla$ of 
the metric and lift it to get the corresponding spin connection $\nabla^s$. 
Referring to the splitting 
${\mathfrak s}{\mathfrak o}(4)={\mathfrak s}{\mathfrak 
u}(2)^+\oplus{\mathfrak s}{\mathfrak u}(2)^-$ the projected connection 
$\nabla_A:=p^+(\nabla^s)$, with $p^+: \so\rightarrow\su^+$ being the 
canonical projection, is an SU$(2)^+$ connection on the SU$(2)^+$ 
rank two complex vector bundle $\Sigma^+$ (the positive chiral spinor 
bundle). 

Assume $(M,g)$ is moreover {\it anti-half-conformally flat} i.e., 
$W^+_g=0$. Picking up an everywhere smooth and positive function $f$ on $M$, 
one can consider the conformally rescaled metric $\tilde{g}:=f^2g$. The new 
space $(M,\tilde{g})$ still satisfies $W^+_{\tilde{g}}=0$ by conformal 
invariance of the Weyl-tensor. Its scalar curvature however 
transforms as $f^3s_{\tilde{g}}=6\triangle_gf+fs_g$ where $\triangle_g$ 
refers to the scalar Laplacian on $(M,g)$. If this new scalar 
curvature $s_{\tilde{g}}$ happens to vanish over $M$ i.e., the 
rescaling function satisfies the equation
\begin{equation}
\triangle_gf+\frac{s_g}{6}f=0
\label{konform-laplace}
\end{equation}
then the corresponding projected connection $\nabla_{\tilde{A}}$, 
constructed like $\nabla_A$ before, satisfies the anti-self-duality 
equations $*_{\tilde{g}}F_{\tilde{A}}=-F_{\tilde{A}}$ over 
$(M,\tilde{g})$ by  a variant of the Atiyah--Hitchin--Singer theorem 
\cite[Proposition 2.2]{ati-hit-sin}. However taking 
into account the conformal invariance of self-duality, 
we obtain that $\nabla_{\tilde{A}}$ in fact satisfies the 
{\it original} SU$(2)^+$ anti-self-duality equations 
$*_gF_{\tilde{A}}=-F_{\tilde{A}}$ over $(M,g)$. Hence we have 
construcred a solution of the anti-self-duality equation via the {\it 
conformal rescaling method}. 

This anti-self-dual connection locally can be constructed as follows. 
Let $U\subset M$ be a coordinate ball and let $(\xi^0,\dots,\xi^3)$ 
denote the orthonormal frame field which diagonalizes $g\vert_U$. In 
this gauge write $\nabla\vert_U =\dd +\omega$ for the Levi--Civita 
connection. The rescaled metric $f^2g$ has a corresponding 
local frame field $(f\xi^0,\dots,f\xi^3)$ and Levi--Civita connection 
$\widetilde{\nabla}\vert_U=\dd +\tilde{\omega}$. 
The connection $1$-form $\tilde{\omega}$ can be calculated by the aid of 
the Cartan equation $\dd (f\xi^i)+\tilde{\omega}^i_j\wedge (f\xi^j)=0$. 
An easy computation shows that 
\[\tilde{\omega}^i_j=\omega^i_j+((\dd\log 
f)\xi_i)\xi^j-((\dd\log f)\xi_j)\xi^i\]
where $(\xi_0,\dots,\xi_3)$ denotes the dual vector field on 
$TU$ (with respect to the original metric $g\vert_U$) and  
$(\dd\log f)\xi_i=g\vert_U({\rm grad} (\log f), \xi_i)$. Lifting 
$\widetilde{\nabla}$ to the spin 
bundle then projecting it to the positive side (using the 't Hooft 
matrices, cf. \cite{ete-hau3}) and exploiting the 
identification $\su^+\cong{\rm Im}\:\HH$ given by 
$(\sigma_1,\sigma_2,\sigma_3)\mapsto (\ii ,-\jj ,-\kk )$ we obtain in 
this induced gauge
\begin{equation}
\tilde{A}=A+\frac{1}{2}{\rm Im}(({\bf d}\log f)\boldsymbol{\xi})
\label{perturbacio}
\end{equation}
with $A:=p^+(\omega^s)$ being the projection of the 
original spin connection $\nabla^s\vert_U=\dd +\omega^s$ and the 
quaternion-valued $0$-form ${\bf d}\log f$ and $1$-form $\boldsymbol{\xi}$ in 
the second term are defined respectively as follows: 
\[{\bf d}\log f:=\dd\log f(-\xi_0+\xi_1\ii +\xi_2\jj +\xi_3\kk ) 
,\:\:\:\:\:\boldsymbol{\xi}:=\xi^0+\xi^1\ii +\xi^2\jj +\xi^3\kk .\]
Writing $a\vert_U:=\frac{1}{2}{\rm Im}(({\bf d}\log 
f)\boldsymbol{\xi})$ for the only term involving $f$, it is clear that it 
extends over $M$ as an ${\rm Im}\:\HH$-valued $1$-form i.e., one finds 
\begin{equation}
a\in C^\infty (M;\Lambda^1M\otimes{\rm End}\:\Sigma^+)
\label{a1}
\end{equation}
and it is insensitive for an overall rescaling of $f$ by a non-zero real 
constant. This implies that from the point of view of our anti-instanton 
factory, we can regard two conformal scaling functions identical if they 
differ by an overall rescaling only. Moreover note that the new 
connection takes the form $\nabla_{\tilde{A}}=\nabla_{A+a}$ i.e., 
this whole conformal rescaling procedure can be regarded as an 
anti-self-dual perturbation of the original connection. 
For an explicit application of this method we refer to 
\cite{jac-noh-reb} over the flat $\R^4$ and to \cite{ete-hau3} in case 
of the multi-Taub--NUT space.

If $(M,g)$ is compact then this machinery provides us with an 
anti-instanton on $\Sigma^+$. If the space is non-compact 
but complete and $\nabla_A$ has finite energy (which can be achieved if the 
metric decays sufficiently fast to the flat one at infinity as we have 
seen) then if the perturbation $a$ also decays fast i.e., $f$ is sufficiently 
bounded at infinity then $\nabla_{\tilde{A}}$ has still finite energy 
consequently we obtain an anti-instanton on the same bundle $\Sigma^+$ 
over $(M,g)$. 

An important particular case is when the original space $(M,g)$ has 
already vanishing scalar curvature. In this case (\ref{konform-laplace}) 
cuts down to the ordinary Laplace equation $\triangle_gf=0$ and the 
original connection $\nabla_A$ is already 
anti-self-dual corresponding to the trivial solution $f=1$.

However it is clear by the maximum principle that in many 
situations there are no non-trivial everywhere smooth bounded solutions to 
(\ref{konform-laplace}) therefore this method is apparently vacuous. 
Fortunately we can modify $f$ to possess mild singularities over $M$. 
At this level of generality by a ``mild singularity'' we mean the following: 
\begin{itemize}
\item[(i)] There are finitely many points $\{y_1,\dots,y_k\}\subset M$ 
such that $f$ diverges in $y_i$ for all $i=1,2,\dots k$;

\item[(ii)] The local energy of the rescaled connection around 
each point is finite i.e., for $U_i$ a neighbourhood of 
$y_i$ one has 
\[\Vert F_{\tilde{A}}\Vert_{L^2(U_i\setminus\{ y_i\})}<+\infty\]
for all $i=1,2,\dots,k$. (In terms of $f$ this 
condition simply controls its behaviour about $y_i$.)
\end{itemize} 
Under these circumstances Uhlenbeck's theorem \cite{uhl} allows us to remove 
the singularities from $\nabla_{\tilde{A}}$ introduced by $f$. Therefore 
if $f$ is a positive solution to (\ref{konform-laplace}) with such 
singularities then the resulting anti-self-dual connection will be smooth 
again and corresponds to a non-trivial anti-instanton over a bundle $E$ if 
$M$ is compact. This SU$(2)^+$ bundle may be no more 
isomorphic to $\Sigma^+$ as a consequence of the allowed singularities 
of the scaling function. In fact its topological type is determined by 
the singularity structure of $f$. In case of a 
non-compact complete base geometry if in addition $f$ is sufficiently 
bounded at infinity then the resulting anti-self-dual connection has 
finite energy i.e., describes an anti-instanton on the same bundle $\Sigma^+$, 
regardless how many singularities $f$ possesses. This is because in 
fact over a non-compact four-manifold all SU$(2)^+$ bundles, including 
$\Sigma^+$, are isomorphic to the trivial one. The significance of the 
singularity structure of scaling functions in instanton theory therefore 
becomes transparent. 

We move on and recall twistor theory of anti-half-conformally flat 
geometries and in particular the Atiyah--Ward construction 
\cite{ati-war} which is another powerful way to solve the 
(anti-)self-duality equations at least in principle.

Let us recall the general theory \cite[Chapter 13]{bes}. Let $(M,g)$ be 
an oriented Riemannian four-manifold and consider the projectivization
$Z:=P(\Sigma^+)$ of the positive chiral spinor bundle $\Sigma^+$ on $M$
(this projectivization is well-defined even if $M$ is not spin).
Clearly, $Z$ admits a fiber bundle structure $p :Z\rightarrow M$ with 
$\C P^1$'s as fibers. A fiber $F_x$ over $x\in M$ can be viewed as the 
parameter space of orthogonal complex structures on $T_xM$ whose induced 
orientations agree with that of $M$. Using the Levi--Civita connection 
on $(M,g)$ we can split the tangent space $T_zZ$ at any point $z\in Z$; 
writing $z=(x,p)$ with $x\in M$ and $p\in F_x$ one endows the horizontal 
subspace $H_{(x,p)}\cong T_xM$ of $T_{(x,p)}Z$ with the complex 
structure $J_p$ given by $p\in F_x$ itself meanwhile the vertical 
subspace $V_{(x,p)}\cong T_pF_x$ of $T_{(x,p)}Z$ is equipped with the 
standard complex structure $I$ coming from $F_x\cong\C P^1$. Therefore 
$T_{(x,p)}Z$ carries a complex structure $(J_p,I)$ and $Z$ as a real 
six-manifold gets possess an almost complex structure. A basic theorem 
of Penrose \cite{pen} or Atiyah, Hitchin and Singer 
\cite[Theorem 4.1]{ati-hit-sin}\cite[Theorem 13.46]{bes}
states that this complex structure is integrable i.e., $Z$ is a complex
three-manifold if and only if $(M,g)$ is anti-half-conformally flat. 
The space $Z$ is called the
{\it twistor space} of $(M,g)$. The complex structure in fact depends
on the confomal class of the metric only. 

The basic holomorphic structure of $Z$ looks like as follows 
\cite[Chapter 13]{bes}:
\begin{itemize}

\item[(i)] For all $y\in M$ the fibers $F_y\subset Z$ 
represent {\it holomorphic lines} $Y\subset Z$ with normal bundles 
isomorphic to $H\oplus H$ where $H$ is the usual line bundle $H$  with 
$\langle c_1(H), [Y]\rangle =1$ on $Y\cong\C P^1$;

\item[(ii)] The antipodal map on $Y$ induces a {\it real structure} on 
$Z$ i.e., an anti-holomorphic involution $\tau : Z\rightarrow Z$ 
satisfying $\tau^2={\rm Id}_Z$;

\item[(iii)] The space of all holomorphic lines in $Z$ near the 
fibers $Y$ form a locally complete family which is a complex 
manifold $M^\C$ of complex dimension four 
and the fibers $Y$ among them are dinstinguished by the property that the real 
structure fixes them, hence are called {\it real lines}.
\end{itemize} 

\noindent A remarkable property of twistor spaces is that they allow one 
to solve certain linear or non-linear field equations over $(M,g)$. 
These equations are of great importance in physics. We shall also quickly 
review these constructions. 

1. The linear field equation relevant to us here is the {\it conformal
scalar Laplace equation} (\ref{konform-laplace}). We shall denote by 
$H_Z$ the holomorphic line bundle on the whole $Z$ such that for any
fiber the restriction $H_Z\vert_Y$ is the standard line bundle $H$ on 
$Y\cong\C P^1$ and by $H^k_Z$ its $k$th tensor product. For notational 
convenience write $\co_Z(k)$ for the sheaf $\co (H^k_Z)$. 
Recall (cf. e.g. \cite{hit2}) that complex-valued real
analytic functions satisfying (\ref{konform-laplace}) on open subsets 
$(U,g\vert_U)\subset (M,g)$ are in one-to-one correspondence with 
elements of the sheaf cohomology group 
$H^1(p^{-1}(U);\co_{p^{-1}(U)}(-2))$ i.e., there is a natural 
isomorphism
\[T:H^1(p^{-1}(U);\co_{p^{-1}(U)}(-2))\cong{\rm
Ker}\:\triangle_g\vert_U\]
called the {\it Penrose transform}: given an element $\varphi\in
H^1(p^{-1}(U);\co_{p^{-1}(U)}(-2))$ and a real line $Y\subset
p^{-1}(U)$ one can take the restriction $\varphi\vert_Y\in
H^1(Y;\co_Y(-2))\cong\C$. This gives rise to a complex-valued
function $f(y):=\varphi\vert_Y$ consequently one has a map
\[T:H^1(p^{-1}(U);\co_{p^{-1}(U)}(-2))\longrightarrow C^\infty (U,\C
).\]
It turns out that $f$ is a solution to (\ref{konform-laplace}) on 
$U\subset M$ and all local solutions arise this way.

2. The non-linear {\it (anti-)self-duality equations} of Yang--Mills 
theory, also can be treated in twistor theory 
(cf. e.g. \cite{ati-hit-sin}): solutions of the SU$(2)$
anti-self-duality equations over $(M,g)$ with approriate orientation can
be converted into certain holomorphic vector bundles on $Z$ as follows.
There is a one-to-one correspondence of gauge equivalence classes of
anti-self-dual SU$(2)$ connections on a rank two
complex vector bundle $E$ on $M$ and isomorphism classes of rank
two complex vector bundles $F$ on $Z$ such that
\begin{itemize}

\item[(i)] $F$ is holomorphic;

\item[(ii)] For any real line $Y\subset Z$ the restricted bundle
$F\vert_Y$ is holomorphically trivial;

\item[(iii)] There exists a conjugate linear map 
$\tilde{\tau}:F\rightarrow F$ lying over the real structure $\tau 
:Z\rightarrow Z$ and satisfies $\tilde{\tau}^2=-{\rm Id}_F$.

\end{itemize}
If $\nabla_A$ is irreducible and $Z$ (that is, $M$) is compact then $F$
is a stable holomorphic bundle on $Z$ and in this case $\tilde{\tau}$ is
unique up to scalar multiplication.

This construction, called the {\it Atiyah--Ward correspondence}, goes
as follows. Take an SU$(2)$ vector bundle $E$ and an anti-self-dual
connection on it over $(M,g)$. Let $F:=p^*E$ be a $C^\infty$ bundle;
by anti-self-duality the curvature of the lifted connection
$p^*\nabla_A$ will be of $(1,1)$-type hence $F$ will be
holomorphic. One checks that properties (ii) and (iii) are also
satisfied. Conversely, take a holomorphic bundle $F$ satisfying (ii) and
(iii). Define $E_y$ to be the space of holomorphic sections over
$F\vert_Y$ (this space is two complex dimensional via (ii)). Property
(iii) induces a symplectic structure on $F$ consequently we can take the
corresponding unique unitary connection $\nabla$ on $F$. By
uniquenss and (ii) this connection is the trivial connection on each
$F\vert_Y$ hence is of the form $p^*\nabla_A$ with a connection
$\nabla_A$ on $E$. Since the curvature of this connection is of
$(1,1)$-type with respect to all complex structures on $T_yM$ for all
$y\in M$ it follows that it is anti-self-dual.

3. Another non-linear equation, namely the {\it (Riemannian) 
anti-self-dual vacuum Einstein equations} of general relativity can also be 
naturally adjusted to twistor theory, known as {\it Penrose' non-linear 
graviton construction} \cite{pen,hit1}. Assume here that $(M,g)$ is 
moreover simply connected and Ricci-flat i.e., simply 
connected and anti-self-dual. Then it 
follows that the induced connection on $\Sigma^+$ is the trivial flat one and 
$Z$ can be retracted onto its particular fiber $F_y$ by parallel 
transport. Therefore topologically $Z\cong M\times S^2$ in this case.
The associated twistor space of $(M,g)$ has in addition to the basic 
holomorphic properties above the following:
\begin{itemize}

\item[(i)] There is a holomorphic fibration $\pi
:Z\rightarrow Y$ over a particular line $Y\cong\C P^1$ with fibers being 
the original space $M$. Consequently $H^k_Z\cong\pi^*H^k$ for all 
$k\in\Z$;

\item[(ii)] Since the complex normal bundle $N_Y$ of each $Y$ is 
isomorphic to $H\oplus H$ it follows that the canonical bundle $K_Z$ 
of $Z$ is isomorphic to $\pi^*H^{-4}$.

\end{itemize}
Conversely, it can be shown \cite[Chapter 13]{bes} that a complex
3-space $Z$ having these properties encodes the conformal class of 
an anti-self-dual simply connected four-space $(M,g)$. 
Indeed, the complexified space $M^\C$ carries a natural conformal structure by
declaring two points $Y', Y''\in M^\C$ to be null-separated if and only
if for the corresponding lines $Y'\cap Y''\not=\emptyset$ in $Z$.
Restricting this conformal structure to the real lines
parameterized by $M\subset M^\C$ one comes up with the confomal class of 
an anti-self-dual structure $(M,g)$. The 
particular metric emerges via an appropriate isomorphism $K_Z\cong\pi^*H^{-4}$.
\vspace{0.1in}

\noindent We have encountered two techniques for constructing 
anti-instantons over an anti-half-con\-for\-mal\-ly flat space. The link 
between them is an observation of Atiyah \cite{ati} which relates solutions to 
(\ref{konform-laplace}) to holomorphic vector bundles on the twistor space. 
We will work out this in detail in Sect. \ref{five}.

In the rest of this section we introduce a 
special class of four-manifolds: suppose from now on that 
$(M,g)$ is a four dimensional {\it asymptotically locally flat} (ALF) 
space as it was defined in \cite{ete-jar}. By a recent powerful theorem of 
Minerbe \cite[Theorem 0.1]{min1} these four dimensional geometries at 
least in the connected, geodesically complete, hyper-K\"ahler case can also 
be characterized by requiring their curvature to be $L^2$ 
with respect to the measure $\frac{r^4}{{\rm Vol}_g(B^4_r(x))}\dd V_g$ 
as well as prescribing their asymptotical volume growth to be $O(r^\nu)$ 
with $3\leq\nu<4$. 

Framed moduli spaces of certain ``admissible'' 
SU$(2)$ (anti-)instantons over these geometries have been investigated 
recently \cite{ete-jar}. By ``admissible'' we mean the following. 
Topologically, an ALF space $(M,g)$ (with a single end) admits a 
decomposition $M=K\cup W$ where $K$ is a compact interior space and $W$ is 
an end or neck homeomorphic to $N\times\R^+$ where $\pi_N :N\rightarrow 
B_\infty$ is a connected, compact, oriented three-manifold fibered over a
compact Riemann surface $B_\infty$ with circle fibers $F\cong S^1$. 
Note that for any $x\in M$ and sufficiently large $R>0$ $W$ is homeomorphic to
$M\setminus B^4_R(x)$ where $B^4_R(x)$ is a geodesic ball about $x\in M$ 
with radius $R>0$.

\begin{definition}
Let $(M,g)$ be an ALF four-manifold. Take an arbitrary smooth, finite 
energy {\rm SU$(2)$} connection $\nabla_A$ on a (necessarily trivial) 
{\rm SU$(2)$} vector bundle $E$ over $M$. This connection is said 
to be {\em admissible} if it satisfies two conditions (cf. \cite{ete-jar}): 
\begin{itemize}
\item[{\rm (i)}] The first is called the {\em weak holonomy 
condition} and says that to $\nabla_A$ there exists a smooth flat 
{\rm SU$(2)$} connection $\nabla_\Gamma\vert_W$ on 
$E\vert_W$ along the end $W\subset 
M$ and a constant $c_1=c_1(g)>0$, independent of $R>0$, such that in any 
smooth gauge on $M\setminus B^4_R(x)$ the inequality 
\[\Vert A-\Gamma\Vert_{L^2_{1,\Gamma}\left( M\setminus 
B^4_R(x)\right)}\leq c_1\Vert F_A\Vert_{L^2\left( M\setminus B^4_R(x)\right)}\]
holds along the neck $M\setminus B^4_R(x)$;

\item[{\rm (ii)}] The second condition requires $\nabla_A$ to {\em decay 
rapidly} at infinity i.e., 
\[\lim\limits_{R\rightarrow\infty}\sqrt{R}\:\Vert 
F_A\Vert_{L^2\left( M\setminus B^4_R(x)\right)}=0.\]
\end{itemize}
\label{feltetelek}
\end{definition}

\noindent The first condition ensures us that the connection has a 
well-defined holonomy at infinity given by $\nabla_\Gamma\vert_W$; 
however note that if the infinity is not simply connected then this 
holonomy can in principle be non-trivial.\footnote{In this paper the 
case of the trivial holonomy will be considered only.} The 
second condition regulates how fast the finite energy connection decays to 
this flat connection and also can be reformulated (cf. \cite{ete-jar}) 
by saying that $F_A$ belongs to a weighted Sobolev space with weight $\delta 
=\frac{1}{2}$ i.e., $F_A\in L^2_{\frac{1}{2}}(M;\Lambda^2M\otimes{\rm 
End}\:E)$. Note that this second condition is stronger than  
the finite energy requirement on $\nabla_A$. It was demonstrated in 
\cite{ete-jar} that irreducible SU$(2)$ anti-instantons satisfying 
conditions (i) and (ii) with a fixed limiting flat connection 
$\nabla_\Gamma\vert_W$ and energy $e$ form framed moduli spaces 
$\cm (e,\Gamma )$ which are smooth finite dimensional manifolds. Their 
dimensions have also been calculated by the aid of a dimension formula.

Now take an ALF space $(M,g)$ with a fixed orientation whose curvature 
satisfies $W_g^+=0$ and $s_g=0$ therefore the Atiyah--Hitchin--Singer 
theorem applies and the projected Levi--Civita connection $\nabla_A$ is 
an anti-self-dual connection on $\Sigma^+$. By an appropriate curvature 
decay imposed on the metric (which is incorporated into the precise 
definition of an ALF space, cf. \cite{ete-jar}) this connection has not 
only finite energy but also decays 
rapidly as in (ii) of Definition \ref{feltetelek}; moreover if 
additionally it satisfies the weak holonomy condition (i) in 
Definition \ref{feltetelek} then it gives rise to an anti-instanton on 
$\Sigma^+$ in the sense above. Performing a conformal rescaling of the 
metric we obtain another anti-self-dual connection 
$\nabla_{\tilde{A}}$ on $\Sigma^+$. Since it is of the 
form $\nabla_{A+a}$ with $a$ as in (\ref{a1}) except finitely many 
points, its curvature reads as $F_{\tilde{A}}=F_{A+a}= F_A+\dd_Aa+a\wedge a$ 
consequently for its energy
\begin{equation}
\Vert F_{\tilde{A}}\Vert_{L^2_{\frac{1}{2}} (M)}\leq \Vert 
F_A \Vert_{L^2_{\frac{1}{2}} 
(M)}+\Vert\dd_Aa\Vert_{L^2_{\frac{1}{2}} (M)}+c_2\Vert 
a\Vert^2_{L^2_{\frac{1}{2},1,\Gamma} (M)}
\label{perturbaltenergia}
\end{equation}
holds with a constant $c_2=c_2(g)>0$. We conclude that the regularity and decay 
properties of the perturbed connection is determined by the perturbation term. 
More precisely, the conditions
\begin{equation}
a\in L^2_{\frac{1}{2},1,\Gamma} (M;\Lambda^1M\otimes{\rm 
End}\:\Sigma^+),\:\:\:\:\:\tilde{A}\in 
L^2_{loc,1,\tilde{A}}(M;\Lambda^1M\otimes{\rm 
End}\:\Sigma^+),\:\:\:\:\:F_{\tilde{A}}\in 
L^2_{loc}(M;\Lambda^-M\otimes{\rm End}\:\Sigma^+)
\label{a2}
\end{equation}
ensure us that $a$ affects neither the asymptotics of $\nabla_A$ 
or its fall-off and any pointlike singularities are removable. Therefore 
if $\nabla_A$ is smooth and obeys both the weak 
holonomy and the rapid decay conditions above then 
$\nabla_{\tilde{A}}=\nabla_{A+a}$ will be also smooth, decays rapidly 
and converges to the same limiting flat connection. More precisely one 
can show that the energy $\tilde{e}$ of $\nabla_{\tilde{A}}$ can differ from 
that of $\nabla_A$ by a non-negative integer only \cite{ete-jar}. 
That is, $[\nabla_A]\in\cm (e, \Gamma )$ implies 
$[\nabla_{\tilde{A}}]\in\cm (\tilde{e}, \Gamma )$ 
and $\tilde{e}-e\in\{ 0,1,2,...\}$. In fact as one expects, 
the difference $\tilde{e}-e$ is governed by the singularity structure of the 
perturbation i.e., of the scaling function $f$ in (\ref{konform-laplace}). 

In the next sections we will apply the conformal rescaling technique and 
twistor theory to understand the moduli space over the multi-Taub--NUT 
geometries of unit energy SU$(2)$ anti-instantons decaying rapidly and 
obeying the weak holonomy condition with respect to the trivial connection 
$\nabla_\Theta$.

 
\section{Review of the multi-Taub--NUT geometry}
\label{three}

In this section we take a closer look of the multi-Taub--NUT spaces 
following \cite{and-kro-leb,leb}. Consider $s>0$ distinct 
points $q_1,\dots,q_s\in\R^3$. We will construct a four-manifold $M_V$ and 
a smooth map $\pi : M_V\rightarrow\R^3$ such that $\pi^{-1}(q_j)$ is a 
point for all $j=1,2,\dots,s$ but $\pi^{-1}(x)\cong S^1$ for all other 
points $x\in\R^3\setminus\{ q_1,\dots,q_s\}$. To begin with, let 
$\pi : U_V\rightarrow \R^3\setminus\{ q_1,\dots,q_s\}$ be 
the principal $S^1$ bundle  whose Chern class is $-1$ when restricted to a 
small sphere of radius $r_j<\min\limits_{k\not=j}\vert q_j-q_k\vert$ around 
$q_j$. Thus $\pi^{-1}(B^3_{r_j}(q_j))$ is diffeomorphic to the punctured 
ball $B^4_j\setminus\{ 0\}\subset\R^4$ in a manner such that the $S^1$ 
action becomes the action of $S^1\subset\C$ on $\C^2\cong\R^4$ by scalar 
multiplication. We then define
\[M_V:=(U_V\sqcup B^4_1\sqcup\dots\sqcup B^4_s)/\sim\]
where $\sim$ means that $B^4_j\setminus\{ 0\}$ is identified with 
$\pi^{-1}(B^3_{r_j}(q_j))$. The map $\pi :U_V\rightarrow\R^3
\setminus\{ q_1,\dots,q_s\}$ extends to a 
smooth map $\pi : M_V\rightarrow\R^3$. Note that there is an $S^1$ 
action on $M_V$ whose fixed points with index $-1$ are exactly the 
$q_j$'s---called {\it NUTs}. If $\ell_{ij}$ is a straight line segment 
joining $q_i$ with $q_j$ then $\pi^{-1}(\ell_{ij})\subset M_V$ is a 
smoothly embedded $2$-sphere whose self-intersection number is $-2$. 
The $2$-spheres $\pi^{-1}(\ell_{j, j+1})$ with $j=1,2,\dots,s-1$ are attached 
together according to the $A_{s-1}$ Dynkin diagram and generate the 
singular cohomology group $H^2(M_V,\Z )\cong\Z^{s-1}$.

Take a real number $c>0$ and let $V:\R^3\setminus\{q_1,\dots,q_s\}
\rightarrow\R^+$ be defined by 
\[V(x):=c +\frac{1}{2}\sum\limits_{j=1}^s\frac{1}{\vert x-q_j\vert}.\]
In this paper we will suppose $c=1$. Then $\triangle V=0$ that is, it is 
a positive harmonic function on $\R^3\setminus\{ q_1,\dots,q_s\}$. The 
cohomology class $[\frac{1}{2\pi}*_3\dd V]\in H^2(U_V,\Z )$ is 
the first Chern class of the bundle $\pi:U_V\rightarrow\R^3\setminus\{ 
q_1,\dots,q_s\}$. Consequently there is a connection $\nabla$ on 
$\pi :U_V\rightarrow \R^3\setminus\{ q_1,\dots,q_s\}$ whose curvature is 
$*_3\dd V$. More precisely, let $\omega\in 
C^\infty(\Lambda^1(\R^3\setminus\{ q_1,\dots,q_s\} ))$ be a real valued 
connection $1$-form on $\R^3\setminus\{ q_1,\dots,q_s\}$ so that 
$*_3\dd V=\frac{1}{\ii}F_\omega =\dd\omega$. The form 
$\omega$ is unique up to gauge transformation since $\R^3\setminus\{ 
q_1,\dots,q_s\}$ is simply connected. If we introduce a coordiante system
$(x^1,x^2,x^3,\tau )$ where $x^i$ are Cartesian coordinates on $\R^3$ 
and $\tau\in [0,2\pi )$ parameterizes the circles on $U_V$ then  
the {\it multi-Taub--NUT metric} $g_V\vert_{U_V}$ on $U_V$ is defined to be 
\begin{equation}
\dd s^2=V\pi^*((\dd x^1)^2+(\dd x^2)^2+(\dd x^3)^2)+
\frac{1}{V}(\dd\tau +\pi^*\omega )^2.
\label{metrika}
\end{equation}
The space $(U_V,g_V\vert_{U_V})$ extends smoothly across the fixed points
of the $S^1$ action: write $V(x)=V_j(x)+f(x)$ with 
$V_j(x):=\frac{1}{2\vert x-q_j\vert}$ and $f$ a smooth function 
around $q_j$. Let $\omega_j$ be the connection whose curvature is  
$*_3\dd V_j$ then we obtain a decomposition
\[\dd s^2=V_j\pi^*((\dd x^1)^2+(\dd x^2)^2+(\dd 
x^3)^2)+\frac{1}{V_j}(\dd\tau +\pi^*\omega_j)^2+h\]
where $h$ is a smooth symmetric bilinear form around $p_j$ induced by 
$f$. Moreover the first term is just the $1$-Eguchi--Hanson metric on
$\R^4$ which is known to be isometric to the standard flat metric on 
$\R^4$ hence extends over $p_j$.

We come therefore up with a Riemannian manifold $(M_V,g_V)$ with an 
isometric $S^1$ action whose fixed points with index $-1$ are $\pi^{-1}(q_j)$. 
It is moreover complete with a single ALF end, is anti-self-dual 
since $\dd\omega =*_3\dd V$ that is, anti-half-conformally flat and 
Ricci-flat. Since $M_V$ is simply connected, it follows that it is 
hyper-K\"ahler as well i.e., there is an entire $2$-sphere's worth of 
complex structures for which $g_V$ is a K\"ahler metric.

To display these complex structures explicitly at least on $U_V$, let 
$e_1,e_2,e_3$ be an orthonormal basis in $\R^3$. Consider these as constant 
vector fields on $\R^3$ and let $\hat{e}_1,\hat{e}_2,\hat{e}_3$ be their 
horizontal lifts to $U_V$ via the connection, as well as let $\hat{e}_0$ 
be the generator of the $S^1$-action. Then 
\begin{equation}
(\xi_0,\xi_1,\xi_2,\xi_3):=\left(\sqrt{V}\hat{e}_0, 
\frac{1}{\sqrt{V}}\hat{e}_1, \frac{1}{\sqrt{V}}\hat{e}_2, 
\frac{1}{\sqrt{V}}\hat{e}_3\right)
\label{bazis}
\end{equation}
is an orthonormal frame on $(U_V,g)$ which we define to be 
oriented. Relative to this frame we take the almost complex structure
\[J_{e_1}:=\begin{pmatrix}
             0 & -1 & 0 & 0\cr
             1 & 0 & 0 & 0\cr
            0 & 0 & 0 & -1\cr
              0 & 0 & 1 & 0
\end{pmatrix}.\]
It is parallel hence is integrable and it depends on the choice of 
$e_1$. Of course it extends over the whole 
$(M_V, g_V)$ which we continue to denote by $J_{e_1}$. By the aid of 
this picture we see that all possible orthogonal complex 
structures with compatible orientation on $(M_V, g_V)$ are parameterized 
by a $\C P^1$ and the situation can be described as a 
holomorphic fibration: the fiber over $e_1\in\C P^1$ is the complex 
manifold $(M_V, J_{e_1})$ in accord with the general theory outlined in 
Sect. \ref{two}. Moreover if 
the direction of any of the straight line segments $\ell_{ij}$ coincides 
with $e_1$ then  the corresponding $2$-sphere $\pi^{-1}(\ell_{ij})$ 
will represent a holomorphic curve with self-intersection number $-2$ in 
$(M_V, J_{e_1})$. For a generic complex structure however, there are no 
holomorphic projective lines. 

These complex structures on the whole $(M_V, g_V)$ can be described in a 
rather explicit algebraic way as follows (cf. \cite{and-kro-leb, leb}). 

\begin{lemma} Fix a generic configuration of points $q_1,\dots,q_s$ in 
$\R^3$ and consider the associated space $(M_V,g_V)$. Take a direction 
$e_1$ in $\R^3$. Assume that $e_1$ is not parallel with any line 
segments $\ell_{ij}$ joining $q_i$ with $q_j$. Then $(M_V, J_{e_1})$ is 
biholomorphic to the algebraic surface $X\subset\C^3$ given by
\begin{equation}
xy-(z-p_1)\dots (z-p_s)=0
\label{model}
\end{equation}
with $(x,y,z)\in\C^3$ and fixed distinct complex numbers $p_1,\dots, 
p_s\in\C$. 

The $S^1$-action on $X$ is given by $(x,y,z,\tau )\mapsto (x{\rm
e}^{\ii\tau}, y{\rm e}^{-\ii\tau}, z)$ with fixed points $(0,0,p_j)\in 
X$ of index $-1$ for all $j=1,\dots,s$. Under this 
biholomorphism the fixed point $\pi^{-1}_V(q_j)\in M_V$ of the 
$S^1$-action on $M_V$ correspond to the fixed point $(0,0,p_j)\in X$ 
with the same index for all $j=1,\dots ,s$ consequently with this 
$S^1$-action on $X$ we obtain an isomorphism of complex manifolds with an 
$S^1$-action.

There exist at most $2\binom{s}{2}=s(s-1)$ different directions when $e_1$ 
is parallel with any $\ell_{ij}$. In this case $(M_V, J_{e_1})$ is 
biholomorphic to the complex surface $\widetilde{X}$ which arises 
by taking $p_i=p_j$ for fixed pair $(i,j)$ in (\ref{model}) and 
blowing up the resulting singular surface $X^*\subset\C^3$ in 
$(0,0, p_j)\in X^*$. 

As a real oriented four-manifold $\widetilde{X}$ is diffeomorphic to $X$.
\label{komplex}
\end{lemma}

\noindent{\it Proof.} By genericity we assume that the straight line 
segments $\ell_{ij}\subset\R^3 $ with $i,j=1,2,\dots, s$ are mutually 
non-parallel. First take a direction $e_1$ in $\R^3$ which is 
not parallel with any $\ell_{ij}$. Let $(e_2e_3)\subset\R^3$ be the 
oriented plane passing through the origin, perpendicular to $e_1$ and  
identify it with the complex plane with its standard orientation: 
$(e_2e_3)\cong\C$ such that $e_2\mapsto 1$ and $e_3\mapsto\ii$. Let 
$p_{e_1}:\R^3\rightarrow\C$ be the projection along $e_1$ 
which induces $P_{e_1}:=p_{e_1}\circ\pi :M_V\rightarrow\C$, the 
projection from $M_V$. Then $P_{e_1}$ is a holomorphic map from $(M_V, 
J_{e_1})$ to $\C$. For a point $a\in M_V$ define a complex coordinate by 
$z:=P_{e_1}(a)\in\C$. In particular for a NUT $q_j\in\R ^3$ put 
$p_j:=P_{e_1}(\pi^{-1}(q_j))=p_{e_1}(q_j)\in\C$ 
for all $j=1,\dots,s$ and also call them NUTs.

Next we analyze the preimage $P^{-1}_{e_1}(z)\subset M_V$. If
$z\not=p_j\in\C$ then $P^{-1}_{e_1}(z)$ is holomorphically isomorphic to
an infinite cylinder $\C^*:=\{ (x,y)\in\C^2\:\vert\: xy=1\}$ in $(M_V,
J_{e_1})$, since the directions $\hat{e}_0$ and $\hat{e}_1$ span a
holomorphic line for the complex structure $J_{e_1}$. On the other hand
$P^{-1}_{e_1}(p_j)$ is homeomorphic to two complex affine lines 
intersecting in one point. In particular we have a model like
$\{ (x,y)\in\C^2\:\vert\: xy=0\}$ for $P^{-1}_{e_1}(p_j)$.
Since the complex structure on the fibers is independent of $z$,
the above two cases can be managed together into a global equation by 
writing $xy=c(z)$ where $c(z)$ is a polynomial in the variable $z$ on 
$\C$ which vanishes exactly in the points $p_1,\dots,p_s\in\C$. 
Therefore we put
\[c(z):=(z-p_1)\dots(z-p_s)\]
to obtain the desired biholomorphism between $(M_V, J_{e_1})$ and $X$
given by (\ref{model}).

The statement about the $S^1$-action is clear. 

If $e_1$ happens to be parallel with any $\ell_{ij}$ connecting 
$q_i\in\R^3$ and $q_j\in\R^3$ then we find 
$p_{e_1}(q_i)=p_i=p_{e_1}(q_j)=p_j\in\C$ hence two roots 
coincide in (\ref{model}). However now 
$P^{-1}_{e_1}(p_i)=P^{-1}_{e_1}(p_j)$ consists of not only two 
complex affine lines as above but in addition the projective line 
$\pi^{-1}(\ell_{ij})$ and the two complex affine lines do not intersect 
each other rather each individual line hits $\pi^{-1}(\ell_{ij})$ in 
one point respectively. These two points are distinct. The situation can be 
modeled on $\{ (x,y,[u:v])\in\C^2\times\C P^1\:\vert\:xy=0, xv=yu\}$. 
Therefore to obtain the resulting 
smooth complex surface we have to take a singular curve $X^*\subset\C^3$ 
represented by (\ref{model}) with $p_i=p_j$ and blow it up at 
$(0,0,p_i)=(0,0,p_j)\in X^*$ to obtain $\widetilde{X}$.
   
It readily follows that as a real oriented four-manifold, $\widetilde{X}$ is 
diffeomorphic to $X$ since both are just the original $M_V$ equipped with 
different complex structures. $\Diamond$
\vspace{0.1in} 

\begin{remark}\rm Incidentally we
note that the above definition of the coordinate $z$
depends on the identification of the oriented plane perpendicular to 
$e_1$ with $\C$, i.e. on the choice of a vector $e_2$ in the tangent 
space of $\C P^1$ at $e_1$. For all such choices, the above procedure yields a
complex number. As $T\C P^1\cong H^2$, this means that globally over
$\C P^1$ the coordinate $z$ must be regarded as a section of the line
bundle $H^2$. This fact will play an important role in constructing 
the twistor space, cf. Sect. \ref{five}.
\end{remark}


\section{Moduli spaces with trivial holonomy}
\label{four}

We already know that a multi-Taub--NUT space $(M_V, g_V)$ is 
anti-half-conformally flat and Ricci-flat. Consequently for its 
curvature $W^+_{g_V}=0$ and $s_{g_V}=0$ holds and the 
Atiyah--Hitchin--Singer theorem applies and we can construct 
SU$(2)^+$ (anti-)self-dual connections on the (trivial) positive chiral 
spinor bundle $\Sigma^+$. If we fix an orientation coming from any compex 
structure in the hyper-K\"ahler family, these connections will be 
anti-self-dual. As we have seen in Section \ref{two}, this approach amounts 
to find positive solutions of the Laplace equation $\triangle_{g_V} f =0$ 
with a finite number of mild singularities and appropriate bounds 
at infinity over the original manifold $(M_V, g_V)$. This program 
partially was carried out in \cite{ete-hau3,ete} and will be 
stretched to its limits here.

An abundance of such harmonic functions are provided by minimal positive 
Green functions therefore we make a short digression on them following 
\cite{gri,heb, var}.

\begin{definition} 
Let $(M,g)$ be a 
non-compact complete Riemannian manifold of dimension $m\geq 2$ and 
$y\in M$ be a point. A function $G(\cdot ,y):M\setminus\{y\}\rightarrow\R$ 
is called the {\em minimal positive Green function concentrated at 
$y\in M$} if it has the following properties:
\begin{itemize}

\item[{\rm (i)}] It satisfies $\triangle_g G(\cdot ,y)=\delta_y$ 
in the sense of distributions on $(M,g)$; 

\item[{\rm (ii)}] Let $r:=d(\cdot ,y)$ denote the distance from $y$ on 
$(M,g)$. If $r\rightarrow +\infty$ then $G(\cdot, y)$ tends to zero and if 
$r\rightarrow 0$ then $G(\cdot, y)$ is $O(r^{2-m})$ and 
$\vert\nabla_y G(\cdot, y)\vert_g$ is 
$O(r^{1-m})$ for $m\geq 3$ as well as $G(\cdot, y)$ is $O(\log r)$ and 
$\vert\nabla_y G(\cdot, y)\vert_g$ is $O(r^{-1})$ for $m=2$;

\item[{\rm (iii)}] Moreover $0<G(\cdot ,y)$ and $G(\cdot ,y)\leq G'(\cdot ,y)$
for any other positive Green function with the same singularity.
\end{itemize}
\label{green}
\end{definition}
\noindent This function if exists is obviously unique and is characterized 
by these properties.

One can try to construct them as follows \cite[p. 229]{heb}. 
Take $y\in M$ and $\Omega\Subset M$ an open subset with smooth boundary 
and $\overline{\Omega}$ compact such that $y\in\Omega$. Let 
$G_\Omega (\cdot ,y)$ be the unique Green function on $\Omega$ satisfying
\[\triangle_gG_\Omega (\cdot ,y)=\delta_y\:\:\:\:\:{\rm
and}\:\:\:\:\:G_\Omega (\cdot ,y)\vert_{\partial\overline{\Omega}}=0.\]
By the maximum principle we can assume that $G_\Omega (\cdot, y)$ is
positive on $\Omega\setminus\{ y\}$ and, if we extend $G_\Omega (\cdot,
y)$ by zero outside $\Omega$, then $G_\Omega (\cdot,
y)\leq G_{\Omega '}(\cdot ,y)$ whenever $\Omega\subseteq\Omega '$.
Consequently setting
\[G(\cdot ,y):=\sup\limits_{\Omega}G_\Omega (\cdot, y)\]
it follows that if this function exists then it is the minimal positive 
Green function with singularity in $y\in M$. One can prove 
\cite[Theorem 8.1]{heb} (also cf. \cite{gri, var}) that minimal positive 
Green functions with singularity $y\in M$ constructed this way either exist 
over $M$ for all $y$ or do not exist at all. 
In the former case $(M,g)$ is called {\it non-parabolic}, while in the 
later case called {\it parabolic}. 

Therefore we should be able to decide whether or not a Riemannian manifold 
is parabolic. For our purposes the easiest way is to recall a result of 
Varopoulos \cite{var} which states that if the Ricci curvature of a 
non-compact complete Riemannian manifold $(M,g)$ with $m\geq 2$ is 
non-negative and for some point $x\in M$ with its geodesic ball of radius $r$
\[\int\limits_1^{+\infty}\frac{r}{{\rm Vol}_g(B^4_r (x))}\:\dd r<+\infty\]  
holds then $(M,g)$ is non-parabolic. 

With this in mind for the Ricci-flat multi-Taub--NUT space $(M_V, g_V)$ we 
find that 
\[{\rm Vol}_{g_V}(B^4_r (x))\sim \frac{8\pi^2}{3}r^3\]
demonstrating that it is non-parabolic (and we can also see that it is 
indeed ALF); consequently for all $y\in M_V$ an associated minimal 
positive Green function $G(\cdot ,y)$ exists. In fact they were explicitly 
constructed by Page \cite{pag} and the simplest of them also appeared in 
\cite{ete-hau3}. Making use of these functions we obtain plenty 
of positive harmonic functions of the form
\begin{equation}
f_{y,\lambda}(x):=1+\lambda G(x,y)
\label{dekompozicio}
\end{equation}
with a real constant $\lambda\in (0,+\infty )$. 

We claim that these functions have the required properties and can be used 
to construct finite-energy smooth SU$(2)^+$ anti-instantons over the 
multi-Taub--NUT space. To check this, we construct the associated 
five parameter family $\nabla_{\tilde{A}_{y,\lambda}}$
of non-gauge equivalent anti-self-dual SU$(2)^+$ connections on
$\Sigma^+$ with respect to the orientation on $(M_V, g_V)$ coming from
the complex structures and demonstrate first that they are 
smooth around $y\in M_V$ hence everywhere and secondly that they are 
admissible hence in particular have finite energy. 

\begin{lemma}
Fix a point $y\in M_V$ and a number $\lambda\in (0,+\infty )$. Take a 
geodesic ball $B^4_\varepsilon (y)\subset M_V$ around $y$ and 
write $B^*_\varepsilon (y):=B^4_\varepsilon (y)\setminus\{ y\}$ for the 
punctured ball. Consider the 
orthonormal frame (\ref{bazis}) on $(B^4_\varepsilon (y), g_V)$. In this
gauge plugging (\ref{dekompozicio}) into (\ref{perturbacio}) write 
$\nabla_{\tilde{A}_{y,\lambda}}\vert_{B^*_\varepsilon (y)}=\dd 
+\tilde{A}_{y,\lambda}$ for the resulting anti-self-dual connection.

There exists an $L^2_{1,\tilde{A}_{y,\lambda}}$ gauge transformation 
$\gamma :B^*_\varepsilon (y)\rightarrow{\rm SU}(2)^+$ such that the 
gauge transformed potential 
$\tilde{A}'_{y,\lambda}:=\gamma^{-1}\tilde{A}_{y,\lambda}\gamma 
+\gamma^{-1}\dd\gamma$ extends smoothly over $B^4_\varepsilon (y)$ that 
is the ${\rm SU}(2)^+$ connection 
$\nabla_{\tilde{A}_{y,\lambda}}\vert_{B^4_\varepsilon (y)}$ 
is smooth on $\Sigma^+\vert_{B^4_\varepsilon (y)}$. 

Consequently the anti-self-dual connetion $\nabla_{\tilde{A}_{y,\lambda}}$ 
is smooth everywhere on $\Sigma^+$. 
\label{mercetrafo}
\end{lemma}

\noindent{\it Proof.} Take an orthonormal frame 
$(\xi^0,\xi^1,\xi^2,\xi^3)$ over $B^4_\varepsilon (y)$ as in 
(\ref{bazis}) with its associated complex structure $J_{e_1}$ on $M_V$. 
Since the multi-Taub--NUT metric is K\"ahler with 
respect to this complex structure, there exists a complex coordinate system 
$(z^1, z^2)$ on $B^4_\varepsilon (y)$ centered in $y$ such that 
\[g_V\vert_{B^4_\varepsilon (y)}=(\delta_{ij}+O(r^2))\dd z^i\dd 
\overline{z}^j\]
i.e., the metric osculates the flat metric in second order. Take a 
neighbourhood $y\in U\subset M_V$ and identify it with $\HH$ such 
that $y\in M_V$ is mapped into ${\bf y}:=y^0+y^1\ii 
+y^2\jj +y^3\kk\in\HH$. Then $B^4_\varepsilon (y)\subset U$ becomes a 
quaternionic ball $B^4_\varepsilon ({\bf y})\subset\HH$ centered about 
${\bf y}$. Introduce a quaternionic coordinate 
${\bf x}\in\HH$ via ${\bf x}:={\bf y}+z^1+z^2\jj$ on it. 
To simplify notation also write $g_V$ for the pullback metric on 
$B^4_\varepsilon ({\bf y})$. There is a bounded positive function $a$ 
such that $r({\bf x})=a({\bf x})\vert{\bf x}-{\bf y}\vert$ hence setting 
$z^1=x^0+x^1\ii$ and $z^2=x^2+x^3\ii$ we find in this real 
coordinate system that
\[\xi^i({\bf x})=(\delta^i_j+O(\vert {\bf x}-{\bf y}\vert^2))\:\dd 
x^j\:\:\:\:\:{\rm and}\:\:\:\:\:\xi_i({\bf x})=(\delta_i^j+O(\vert 
{\bf x}-{\bf y}\vert^2)\frac{\partial}{\partial x^j}.\]
Similarly for the harmonic function in (\ref{dekompozicio}) we obtain
\[f_{{\bf y},\lambda}({\bf x})=1+\frac{\lambda}{\vert 
{\bf x}-{\bf y}\vert^2}+O(\vert {\bf x}-{\bf y}\vert^{-1})\:\:\:\:\:{\rm 
and}\:\:\:\:\:\dd f_{{\bf y} ,\lambda}({\bf x})=-2\lambda\left(\frac{x^i-y^i}
{\vert {\bf x}-{\bf y}\vert^4}\delta^i_j+O(\vert {\bf x}-{\bf 
y}\vert^{-2})\right)\dd x^j.\]
Inserting these expansions into (\ref{a1}) we obtain 
$a_{{\bf y},\lambda}=B_{{\bf y},\lambda}+b_{{\bf y},\lambda}$ for the 
perturbation where the singular Euclidean term looks like
\[B_{{\bf y},\lambda}({\bf x})={\rm
Im}\frac{\lambda (\overline{{\bf x}}-\overline{{\bf y}})\:\dd
{\bf x}}{\vert {\bf x}-{\bf y}\vert^2(\lambda +\vert
{\bf x}-{\bf y}\vert^2)}\]
and $b_{{\bf y},\lambda}$ is of $O(1)$ hence regular in 
$B^4_\varepsilon ({\bf y})$.

Regarding the spin connection $\nabla_A$, note that the SU$(2)^+$ bundle 
$\Sigma^+$ is trivial over $M_V$. The metric is anti-self-dual 
hence the spin connection is flat and since 
$M_V$ is simply connected it is just the trivial flat connection. 
Writing $\nabla_A\vert_{B^4_\varepsilon (y)} =\dd +A$ we simply 
find $A=0$ in the natural gauge (\ref{bazis}). 

Putting all of these into (\ref{perturbacio}) we 
eventually come up in the gauge (\ref{bazis}) we use with an expansion 
of the vector potential as follows:
\[\tilde{A}_{{\bf y},\lambda}=B_{{\bf y},\lambda}+b_{{\bf 
y},\lambda}\]
whose only singular term is $B_{{\bf y},\lambda}$ but 
still $\tilde{A}_{{\bf y},\lambda}\in 
L^2_{1,\tilde{A}_{{\bf y},\lambda}}(\Lambda^1B^*_\varepsilon 
({\bf y})\otimes{\rm End}\Sigma^+\vert_{B^*_\varepsilon ({\bf y})})$. 

Calculating the curvature we find
\[F_{\tilde{A}_{{\bf y},\lambda}}=F_{B_{{\bf y},\lambda}}+\dd_{B_{{\bf 
y},\lambda}}b_{{\bf y},\lambda}+b_{{\bf y},\lambda}\wedge b_{{\bf y},\lambda}\]
that is, 
\begin{equation}
F_{\tilde{A}_{{\bf y},\lambda}}({\bf 
x})=\frac{\lambda\:\dd\overline{{\bf 
x}}\wedge\dd{\bf x}}{(\lambda +\vert{\bf x}-{\bf y}\vert^2)^2}
+\dd_{B_{{\bf y},\lambda}({\bf x})}b_{{\bf 
y},\lambda}({\bf x})+b_{{\bf 
y},\lambda}({\bf x})\wedge b_{{\bf y},\lambda}({\bf x})
\label{gorbulet}
\end{equation}
demonstrating that the curvature $F_{B_{{\bf
y},\lambda}}$ of the Euclidean term as well as $b_{{\bf 
y},\lambda}\wedge b_{{\bf y},\lambda}$ are smooth in the ball 
$B^4_\varepsilon ({\bf y})$ 
meanwhile $\dd_{B_{{\bf y},\lambda}}b_{{\bf y},\lambda}$ is 
singular. Let $\langle\cdot ,\cdot\rangle$ denote the pointwise scalar 
product and $\vert\cdot\vert$ the induced pointwise norm on $\su^+$-valued 
2-forms associated to the Killing norm on $\su^+$ and the metric on
$\Lambda^2M_V$. This gives rise to the
scalar product $(\cdot,\cdot )_{L^2}$ and norm $\Vert\cdot\Vert_{L^2}$ 
on the punctured ball. Separating the regular and singular terms we obtain
\begin{eqnarray}
\Vert F_{\tilde{A}_{{\bf y},\lambda}}\Vert^2_{L^2\left( B^*_\varepsilon 
({\bf y})\right)} & =& \Vert F_{B_{{\bf y},\lambda}} +b_{{\bf 
y},\lambda}\wedge b_{{\bf y},\lambda}\Vert^2_{L^2\left( B^*_\varepsilon 
({\bf y})\right)}+\Vert \dd_{B_{{\bf y},\lambda}}b_{{\bf y},\lambda}
\Vert^2_{L^2\left( B^*_\varepsilon ({\bf y})\right)}\nonumber\\
         & &  \nonumber\\
&+& 2 \left( F_{B_{{\bf y},\lambda}} +b_{{\bf y},\lambda}\wedge b_{{\bf
y},\lambda}\:,\:\dd_{B_{{\bf y},\lambda}}b_{{\bf
y},\lambda}\right)_{L^2\left( B^*_\varepsilon ({\bf 
y})\right)}.\nonumber
\end{eqnarray}
Since the dual of the volume form $*(\dd
V_{g_V}\vert_{B_\varepsilon ({\bf y})}({\bf x}))$ with respect to the 
flat Hodge star $*$ on $\HH$ is $O(\vert{\bf x}-{\bf y}\vert^3)$ 
we find on the one hand that the regular energy term is of course 
finite:
\[\Vert F_{B_{{\bf y},\lambda}} +b_{{\bf y},\lambda}\wedge b_{{\bf 
y},\lambda}\Vert^2_{L^2\left( B^*_\varepsilon ({\bf y})\right)}
=\int\limits_{B^*_\varepsilon ({\bf y})}\vert F_{B_{{\bf 
y},\lambda}}+b_{{\bf y},\lambda}\wedge b_{{\bf
y},\lambda}\vert^2\:\dd V_{g_V}\leq 
c_3\int\limits_0^\varepsilon r^3\dd r <+\infty .\]
The singular term $\vert \dd_{B_{{\bf y},\lambda} ({\bf x})}
b_{{\bf y},\lambda}({\bf x})\vert$ on the other hand is 
$O(\vert{\bf x}-{\bf y}\vert^{-1})$ consequently
\[\Vert \dd_{B_{{\bf y},\lambda}}b_{{\bf
y},\lambda}\Vert^2_{L^2\left( B^*_\varepsilon ({\bf y})\right)}
=\int\limits_{B^*_\varepsilon ({\bf y})}\vert \dd_{B_{{\bf 
y},\lambda}}b_{{\bf y},\lambda}\vert^2\:\dd V_{g_V}\leq
c_4\int\limits_0^\varepsilon r^{-2}r^3\dd r <+\infty .\]
Finally, the cross term also satisfies the estimate (here 
and only here $\vert\cdot\vert$ is the ordinary pointwise absolute value 
on reals)
\begin{eqnarray}
\left\vert\left( F_{B_{{\bf y},\lambda}} +b_{{\bf
y},\lambda}\wedge b_{{\bf
y},\lambda}\:,\:\dd_{B_{{\bf y},\lambda}}b_{{\bf
y},\lambda}\right)_{L^2\left( B^*_\varepsilon ({\bf y})\right)}\right\vert
&\leq &\int\limits_{B^*_\varepsilon ({\bf y})}\left\vert\left\langle 
F_{B_{{\bf y},\lambda}} +b_{{\bf y},\lambda}\wedge b_{{\bf
y},\lambda}\:,\:\dd_{B_{{\bf y},\lambda}}b_{{\bf
y},\lambda}\right\rangle\right\vert\dd V_{g_V}\nonumber\\
&\leq &c_5\int\limits_0^\varepsilon r^{-1}r^3\dd r <+\infty .\nonumber
\end{eqnarray}
We conclude that the last two conditions in (\ref{a2}) hold hence the 
apparent singularity in 
$\nabla_{\tilde{A}_{y,\lambda}}$ at $y$ can be removed by a gauge 
transformation via Uhlenbeck's theorem \cite{uhl} providing us with a 
smooth anti-self-dual ${\rm SU}(2)^+$ connection on the whole positive 
spinor bundle $\Sigma^+$. $\Diamond$
\vspace{0.1in}

\noindent Proceeding further we demonstrate that our solutions are 
anti-instantons, i.e. have finite energy. To carry out this we 
construct an asymptotic expansion for (\ref{dekompozicio}) to see explicitly 
its fall-off. Intuitively \cite{ete}, a harmonic function with a fixed 
singularity on $(M_V, g_V)$ asymptotically looks like pulling the singular 
point of the function as well as all the NUTs in $M_V$ together. Therefore our 
strategy will be as follows: we construct an expansion of a harmonic 
function whose singularity coincides with that of a singular 
multi-Taub--NUT space. After finding its asymptotic expansion we show 
that the decay rate of the leading term remains unchanged if we 
perturb this singular space into the original smooth one. 

In light of Lemma \ref{komplex} this collapsed space arises when all the 
roots of the right hand side in (\ref{model}) coincide: $p_1=\dots =p_s=0$ 
hence the resulting space $X^*\subset\C^3$ is given by $xy-z^s=0$. 
On the other hand, if $(u,v)\in\C^2$ then the polynomials $u^s, v^s$ and 
$uv$ are invariant under the standard action of $\Z_s\subset{\rm SU}(2)$ 
on $\C^2$ and putting $x:=u^s$, $y:=v^s$, $z:=uv$ then they also satisfy 
$xy-z^s=0$ i.e., $X^*$ is a model in $\C^3$ for the $A_{s-1}$-singularity 
$M^*_V:=\C^2/\Z_s$.

Removing the origin $\{ 0\}\in\C^2$ the resulting 
smooth space $W^*:=\C^2\setminus\{ 0\}/\Z_s$ is the neck at 
infinity and is topologically a lens space: $W^*\cong L(s,-1)\times\R^+$. 
Here $L(s,-1)\cong S^3/\Z_s$ denotes the usual lens space of 
$(s,-1)$-type with its orientation inherited from the complex structure on 
$W^*$. Regarding $W^*$ as a complex line bundle of Chern class $-s$ over 
$S^2$ let us denote by $(r,\tau )$ polar coordinates on the fibres 
centered in the singular origin of $\C^2/\Z_s$ and by $(\Theta ,\phi)$ 
the usual spherical coordinates on the base.\footnote{In other words 
$(\tau ,\Theta ,\phi )$ are Euler coordinates on the lens space and 
$r$ is the radial coordiante on the cone $\C^2/\Z_s$. For clarity we 
mention that $r$ defined this way asymptotically 
accords with the distance function $r$ used in Definition \ref{green} 
while $\tau$ here is $\frac{2}{s}$-times the periodic coordinate $\tau$ 
appeared in (\ref{metrika}).} We simply get
\[V(r)=1+\frac{s}{2r},\:\:\:\:\:\alpha (\Theta, \phi)=
\frac{s}{2}\cos\Theta\dd\phi\]
consequently the multi-Taub--NUT metric (\ref{metrika}) 
reduces to a singular metric $g^*_V$ on $M^*_V$ whose shape on $W^*$ is
\[\dd s^2=\frac{2r+s}{2r}(\dd 
r^2+r^2\dd\Theta^2+r^2\sin^2\Theta\dd\phi^2)+
\frac{s^2}{4}\frac{2r}{2r+s}(\dd\tau +\cos\Theta\dd\phi )^2\]
with parameters
\[0<r<+\infty ,\:\:\:\:\:0\leq\tau <\frac{4\pi}{s} ,\:\:\:\:\:0\leq\phi 
<2\pi ,\:\:\:\:\:0\leq \Theta <\pi .\]
One easily calculates 
$\det g^*_V=\frac{s^2}{4}(\frac{2r+s}{2})^2r^2\sin^2\Theta$ 
and by referring to the local expression for the Laplacian of any metric
\[\triangle_g =\sum\limits_{i,j}\frac{1}{\sqrt{\det g}}
\frac{\partial}{\partial x^i}\left(\sqrt{\det g}\:g^{ij}
\frac{\partial}{\partial x^j}\right) ,\]
where $g^{ij}$ are the components of the inverse matrix, our singular
Laplacian looks like
\begin{eqnarray}
\triangle_{g^*_V} & = & \frac{2r}{2r+s}\frac{\partial^2}{\partial r^2}
+\frac{4}{2r+s}\frac{\partial}{\partial r}+\frac{2(2r+s)^2\sin^2\Theta 
+2s^2\cos^2\Theta}{rs^2 
(2r+s)\sin^2\Theta}\frac{\partial^2}{\partial\tau^2}\nonumber\\
&+&\frac{2}{r(2r+s)}\left(\frac{\partial^2}{\partial\Theta^2}+
\cot\Theta\frac{\partial}{\partial\Theta}\right)+\frac{2}
{r(2r+s)\sin^2\Theta}\left(\frac{\partial^2}{\partial\phi^2}
-2\cos\Theta\frac{\partial^2}{\partial\tau\partial\phi}\right) .\nonumber
\end{eqnarray}
This apparently new Laplacian can be re-expressed in old terms as 
follows. The lens space with its standard metric induced from the round 
$S^3$ has an associated Laplacian
\[\triangle_{L(s,-1)} 
=\frac{\partial^2}{\partial\Theta^2}+\cot\Theta\frac{\partial}
{\partial\Theta}+\frac{1}{\sin^2\Theta}\left(\frac{\partial^2}
{\partial\phi^2}-2\cos\Theta\frac{\partial^2}{\partial\tau\partial\phi}+
\frac{\partial^2}{\partial\tau^2}\right)\]
in Euler coordinates. By the aid of this we plainly obtain 
\[V(r)\triangle_{g^*_V}=\frac{\partial^2}{\partial
r^2}+\frac{2}{r}\frac{\partial}{\partial r}+\left(\frac{4}{rs}+ 
\frac{4}{s^2}\right)\frac{\partial^2}{\partial\tau^2}+
\frac{1}{r^2}\triangle_{L(s,-1)}.\]
As we have seen, we have to find positive solutions to 
$\triangle_{g^*_V}f^*=0$ or equivalently to
\begin{equation}
V\triangle_{g^*_V}f^*=0
\label{egyenlet}
\end{equation} 
with singularity in the origin $r=0$. These functions can be constructed 
easily by separating variables. To this end we introduce an orthonormal 
system of bounded smooth functions on the lens space $L(s,-1)$ by the 
(three dimensional) spherical harmonics $Y^{k,l}_j$ with $j=0,1,\dots$ and 
$-j\leq k\leq j$, $-[\frac{2j}{s}]\leq l\leq [\frac{2j}{s}]$ (for a 
fixed $j$, there are $(2j+1)([\frac{2j}{s}]+1)$ independent spherical 
harmonics on $L(s,-1)$). In Euler coordinates these take the shape
\[Y^{k,l}_j(\tau, \phi ,\Theta )=N^{k,l,s}_j\:{\rm 
e}^{\ii\left(\frac{ls}{2}\tau +k\phi\right)} P^{k,l}_j(\cos\Theta )\]
with appropriate normalizing constants $N^{k,l,s}_j$. The 
generalized associated Legendre polynomial $P^{k,l}_j(x)$ is the 
unique real analytic solution on $[-1,1]$ of the following differential 
equation:
\[\left[ (1-x^2)\frac{\dd^2}{\dd x^2}-2x\frac{\dd}{\dd 
x}-\frac{k^2-x(kls)+\frac{l^2s^2}{4}}{1-x^2}+j(j+1)\right] 
P^{k,l}_j(x)=0.\]
Note that taking $l=0$ the function $P^{k,0}_j(x)$ 
reduces to the standard associated Legendre polynomial $P^k_j(x)$ as 
well as $Y^{k,0}_j(\tau ,\phi ,\Theta )$ gives back the usual spherical 
harmonics $Y^k_j(\phi ,\Theta )$ on $S^2$. The spherical harmonics obey
\[\frac{\partial Y^{k,l}_j}{\partial\tau} =\frac{\ii ls}{2}\:Y^{k,l}_j
\:\:\:\:\:{\rm 
and}\:\:\:\:\:\triangle_{L(s,-1)}Y^{k,l}_j=-j(j+1)Y^{k,l}_j.\]
In this basis we seek solutions in the form 
\[f^*(r,\tau ,\phi ,\Theta) 
=\sum\limits_{j=0}^\infty\sum\limits_{k=-j}^j\sum
\limits_{l=-[\frac{2j}{s}]}^{[\frac{2j}{s}]}\lambda^{k,l}_j\varrho^{k,l}_j(r)
Y^{k,l}_j(\tau, \phi ,\Theta )\]
where $\lambda^{k,l}_j$ are complex numbers and $\varrho^{k,l}_j$ are 
radial functions. Inserting this ansatz into (\ref{egyenlet}) and then 
exploiting the eigenfunction properties of the spherical harmonics our 
equation cuts down to a family of ordinary differential equations 
\begin{equation}
\left[\frac{\dd^2}{\dd r^2}+\frac{2}{r}\frac{\dd}{\dd r}-
\frac{j(j+1)}{r^2}-\frac{sl^2}{r}-l^2\right]\varrho^{k,l}_j(r)=0
\label{radial2}
\end{equation}
with $j=0,1,\dots$ and $-j\leq k\leq j$, $-[\frac{2j}{s}]\leq l\leq 
[\frac{2j}{s}]$. The functions $\varrho^{k,l}_j$ are in fact independent 
of the index $k$ due to the U$(2)$ symmetry of the space $(M^*_V,g^*_V)$. 
Two particular solutions if $l=0$ are  
\[K^{0}_{-j-1}(r)=r^{-j-1}\:\:\:\:\:{\rm and}\:\:\:\:\:K^{0}_j(r)=r^j\]
moreover for $l>0$ and $r<\!\!<1$ we know that
\[K^{l}_{-j-1}(r)=r^{-j-1}(1+O(r))\:\:\:\:\:{\rm and}
\:\:\:\:\:K^{l}_j(r)=r^j(1+O(r))\]
meanwhile for $l>0$ and $r>\!\!>1$ we obtain
\[K^{l}_{-j-1}(r)={\rm e}^{-lr}r^{-\frac{sl}{2}-1}(1+O(r^{-1}))
\:\:\:\:\:{\rm and}\:\:\:\:\:K^{l}_j(r)={\rm e}^{lr}r^{\frac{sl}{2}-1}
(1+O(r^{-1})).\]
Separating variables therefore forces us to put the 
pointlike singularity of $f^*$ into the origin of $M^*_V$ showing our 
approach is consistent. To summarize, the general solution formally looks like
\begin{equation}
f^*(r,\tau ,\phi ,\Theta )=\sum\limits_{j=0}^\infty\sum\limits_{k=-j}^j\sum
\limits_{l=-[\frac{2j}{s}]}^{[\frac{2j}{s}]}\left(\lambda^{k,l}_j 
K^{l}_{-j-1}(r)+\mu^{k,l}_j K^{l}_j(r)\right) Y^{k,l}_j(\tau, \phi ,\Theta )
\label{eloallitas}
\end{equation}
with the only singularity in the origin. 

Because $\overline{Y^{k,l}_j}=Y^{-k, -l}_j$, a real basis for the solutions 
is given by 
\[\frac{1}{2}(Y^{k,l}_j+Y^{-k,-l}_j)\:\:\:\:\:{\rm and}\:\:\:\:\:
\frac{1}{2\ii}(Y^{k,l}_j-Y^{-k,-l}_j).\]
Moreover from (\ref{radial2}) we obtain that 
$K_{-j-1}^{l}=K^{-l}_{-j-1}$ and $K^{l}_j=K^{-l}_j$ consequently 
if we suppose $\overline{\lambda^{k,l}_j}=\lambda_j^{-k, -l}$ as well as 
$\overline{\mu^{k,l}_j}=\mu_j^{-k, -l}$ then we come up with formal real 
solutions in (\ref{eloallitas}). In addition asymptotically bounded 
solutions from (\ref{eloallitas}) arise by setting $\mu^{k,l}_j=0$ for $j>0$.

Taking into account that $Y^{0,0}_0$ is a 
constant, the leading $j=0$ term shows that the asymptotic shape of 
asymptotically bounded real harmonic functions looks like $f^*(r)=\mu +
\lambda r^{-1}+O(r^{-2})$ (we put simply $\mu^{0,0}_0=\mu$ and 
$\lambda^{0,0}_0=\lambda$, both real) on the singular space $(M^*_V, g^*_V)$. 

Now we ask ourselves how this picture changes over the original 
non-singular space $(M_V, g_V)$ for the Green function $G(\cdot ,y)$ 
with singularity at $y\in M_V$. For this purpose, we use results of 
\cite{min2}, where a careful study of analytic properties of the Laplace 
operator on perturbations of ALF spaces has been carried out. 
If $W^*$ and $W$ are the necks of $M^*_V$ and $M_V$ respectively then 
we suppose $W^*\cong W$ and use the same coordinate system along them. 
It is clear that for $\alpha =2$ we have 
$g_V\vert_W=g^*_V\vert_{W^*}+O(r^{-\alpha})h$ with some $O(1)$ 
symmetric tensor field $h$. Let $\chi$ be a cut-off function supported 
in a compact neighborhood of $y$ and identically equal to $1$ in a smaller 
neighborhood of $y$. We will try to look for the Green function in the form 
\[G(\cdot ,y) = \chi G^{(4)}(\cdot ,y)+ u\]
where $G^{(4)}(\cdot ,y)$ stands for the Green function on $\R^4$ and 
$u$ is some locally $L^2$ function. Then, $u$ must satisfy the equation 
\begin{equation}\label{inhomLapl}
   \triangle_{g_V} u = f 
\end{equation}
for some prescribed compactly supported smooth function $f$ only depending 
on $\chi$ and $G^{(4)}(\cdot ,y)$. In particular, $f$ belongs to any 
weighted Sobolev space $L^2_\delta(M_V;\R )$. Setting $\delta=2$ it 
follows from \cite[Corollary 2]{min2} that equation (\ref{inhomLapl}) 
can be solved with $u\in L^2_2(M_V;\R )$. Put now 
$\delta':=-\frac{1}{2}+\varepsilon$ for any $\varepsilon>0$, and denote the 
interval $(\delta',\delta)$ by $I$. Notice that the only critical weights in 
$I$ are $\delta_0=\frac{3}{2}$ and $2-\delta_0=\frac{1}{2}$. 
It then follows from \cite[Proposition 4 and Corollary 1]{min2} that up to 
perturbations in $L^2_{-\frac{1}{2}+\varepsilon}(M_V;\R )$ and 
$L^2_{-\frac{3}{2}+\varepsilon ',2}(M_V;\R )$ with some 
$\varepsilon '>0$, $u$ is equal to a linear combination of solutions 
\[K^0_0(r)=r^{\nu_0^+} = 1\:\:\:\:\:{\rm and}\:\:\:\:\: 
K^0_{-1}(r)=r^{\nu_0^-} = r^{-1}\]
of the non-perturbed Laplacian. Notice that the corresponding spherical 
harmonic $Y^{0,0}_0=\frac{\sqrt{s}}{4\pi}$ is constant and 
that the decay corresponding to $L^2_{-\frac{1}{2}+\varepsilon}(M_V;\R )$ 
is $r^{\nu}$ with $\nu = -2 + \varepsilon$, moreover the decay corresponding 
to $L^2_{-\frac{3}{2}+\varepsilon ',2}(M_V;\R )$ is even stronger. 
Consequently (\ref{dekompozicio}) asymptotically also looks like
\begin{equation}
f_{y,\lambda}(r)=1+\lambda r^{-1}+O(r^{-2 + \varepsilon})
\label{alapfv}
\end{equation}
and the leading term is of course positive for $\lambda >0$. One can 
also immediately read off this from the explicit Green function of Page 
\cite[Equation (18)]{pag}. 

By the aid of these observations we obtain

\begin{lemma} The smooth anti-self-dual
connection $\nabla_{\tilde{A}_{y,\lambda}}$ with $y\in M_V$ and 
$\lambda\in (0,+\infty )$ is subject to both the weak 
holonomy condition with respect to the trivial flat 
connection $\nabla_\Theta$ and the rapid 
decay condition of Definition \ref{feltetelek} consequently it has 
not only finite energy but is moreover admissible i.e., its gauge 
equivalence class gives rise to an admissible anti-instanton on 
$\Sigma^+$. The energy of this anti-instanton is of one unit. 
\label{aszimptotikus}
\end{lemma}
\noindent {\it Proof.} We already have seen that the spin connection on 
$\Sigma^+$ is gauge equivalent to the trivial flat connection 
$\nabla_\Theta$.

Consider the neck $W\subset M_V$ of the multi-Taub--NUT
manifold and using the gauge (\ref{bazis}) along $W$
write $\nabla_{\tilde{A}_{y,\lambda}}\vert_W=\dd
+\tilde{A}_{y,\lambda}$ with $\tilde{A}_{y,\lambda}=a_{y,\lambda}$ as 
in (\ref{perturbacio}) where $a_{y,\lambda}$ is the perturbation term 
involving (\ref{dekompozicio}). The asymptotic shape of 
(\ref{dekompozicio}) looks like (\ref{alapfv}) in light of our previous 
calculations. This ensures us that if $\lambda$ 
is finite and positive, $\vert\dd\log f_\lambda\vert_{g_V}$ as well as 
$\vert\nabla_{\Theta}(\dd\log f_\lambda )\vert_{g_V}$ are $O(r^{-2})$ 
along $W$. Consequently, since 
\[\vert a_{y,\lambda}\vert 
=\sqrt{\frac{3}{2}}\:\vert\dd\log f_{y,\lambda}\vert_{g_V}\]
in the gauge (\ref{bazis}), we also obtain that both $\vert 
a_{y,\lambda}\vert$ and $\vert\dd_{\Theta}a_{y,\lambda}\vert$ are 
$O(r^{-2})$ along $W$ for all $y\in M_V$ and $0<\lambda <+\infty$. We 
see that the first condition in (\ref{a2}) is also satisfied. 
Inserting this and $F_\Theta =0$ into (\ref{perturbaltenergia}) we can see 
that our connection can be viewed as a rapidly decaying perturbation of the 
trivial connection hence satisfies both the weak holonomy condition with 
respect to the trivial flat connection as well as the rapid decay condition.

It follows \cite[Theorem 4.2]{ete-jar} that in this situation the 
energy of the connection $\nabla_{\tilde{A}_{y,\lambda}}$ is 
integer hence it must be one because the harmonic function in 
(\ref{dekompozicio}) has only one singularity. $\Diamond$
\vspace{0.1in}

\begin{remark}\rm Before summarizing our findings we may raise the 
question: what about the 
higher overtones with $j>0$ in the harmonic expansion 
(\ref{eloallitas})? The radial function $r^2=|u|^2+|v|^2$ in $\C^2$ 
is invariant under $\Z_s$ consequently we can talk about 
lens spaces of radius $r$ centered at the origin of $\C^2/\Z_s$. Take 
such a lense space with its metric $g$ inherited from the round metric 
of $S^3$. By orthogonality the spherical harmonics obey for $j>0$ that
\[\frac{1}{2}\int\limits_{L(s,-1)}\left( Y^{k,l}_j+Y^{-k,-l}_j\right)\dd
V_g=\frac{1}{2\ii}\int\limits_{L(s,-1)}\left(Y^{k,l}_j-Y^{-k,-l}_j\right)
\dd V_g=0\]
hence these real functions must change sign somewhere
along any lense space showing that in fact all individual terms in 
(\ref{eloallitas}) with $j>0$ vanish on three dimensional subsets of 
$M^*_V$. Additionally we know from the asymptotic formulae for the 
radial functions $K^l_{-j-1}$ and $K^l_j$ that they behave very differently for 
small or large $r$'s consequently the aforementioned 
three dimensional zero sets of the individual higher overtones cannot be 
removed by forming particular linear combinations of them. 
Consequently these harmonic functions represent ``non-mild'' 
singularities hence cannot be used to construct anti-instantons.
\end{remark}
 
\noindent Let $\nabla_\Gamma\vert_W$ be a smooth flat 
${\rm SU(2)}^+$ connection on $\Sigma^+\vert_W$. Then following 
\cite{ete-jar} if $C_0^\infty (M_V;{\rm End}\Sigma^+)$ is the space of 
smooth, compactly supported endomorphisms of the positive chiral spin bundle 
then define the $L^2$ gauge group $\cg_{\Sigma^+}$ as the 
$L^2_{2,\Gamma}$ completion of the space 
\[\{\gamma -1\in C_0^\infty (M_V;{\rm
End}\Sigma^+)\:\vert\:\Vert\gamma 
-1\Vert_{L^2_{2,\Gamma}\left( M_V\right)}<+\infty 
,\gamma\in\mbox{$C^\infty (M_V;{\rm Aut}\Sigma^+)$ a.e.}\}\]
Let us denote by $\cm (e,\Gamma )$ the {\it framed moduli space}
consisting of pairs $([\nabla_A],\Gamma )$, the $\cg_{\Sigma^+}$ 
gauge equivalence classes of smooth irreducible ${\rm SU(2)}^+$ 
anti-self-dual connections on $\Sigma^+$ of energy $e<+\infty$, decaying 
rapidly and obeying the weak 
holonomy condition with respect to $\nabla_\Gamma\vert_W$ and a fixed 
gauge $\Gamma$ at infinity preserved by $\cg_{\Sigma^+}$. 
Associated to this space one has the {\it unframed moduli space} 
$\widehat{\cm}(e, \Gamma )$ consisting of the $\cg_{\Sigma^+}$ 
equivalence classes $[\nabla_A]$ only. In other words it is formed by 
dividing via the group 
\begin{equation}
\widehat{\cg}_{\Sigma^+}\cong\cg_{\Sigma^+}\times{\rm SU}(2)^+
\label{mercecsoport}
\end{equation}
of gauge transformations tending to any fixed element 
$\gamma_0\in{\rm SU}(2)^+$ (not to the identity only) at infinity.

In particular we can consider the spaces $\cm (k,\Theta )$ and 
$\widehat{\cm}(k, \Theta )$ corresponding to the 
trivial flat connection $\nabla_\Theta$; in this case the energy is an 
integer $k\in\N$, cf. \cite[Theorem 4.2]{ete-jar}. 
It was demostrated in \cite[Theorem 3.2]{ete-jar} that these spaces are 
smooth manifolds and \cite[Theorem 4.2]{ete-jar} ensures us that 
$\dim\cm (k,\Theta )=8k$ and therefore $\dim\widehat{\cm}(k,\Theta 
)=8k-3$.
 
Putting together our findings sofar as well as noting that all of our
anti-instantons are irreducible \cite[Theorem 5.1]{ete-hau3} we obtain
the following

\begin{theorem}
Let $(M_V, g_V)$ be the multi-Taub--NUT space with $s>0$ NUTs equipped   
with an orientation coming from any of its complex strucures in the
hyper-K\"ahler family.

Then there exists a $5$ real parameter family of unframed, non-gauge 
equivalent, irreducible, smooth ${\rm SU}(2)^+$ anti-instantons
\[\{ [\nabla_{\tilde{A}_{y,\lambda}}]\:\vert\:(y,\lambda )\in
M_V\times (0,+\infty )\}\]
on the positive chiral spinor bundle $\Sigma^+$ which satisfy the
weak holonomy condition with respect to
the trivial flat connection and also satisfy the rapid decay condition. 
All these anti-instantons have unit energy and provide us with an open 
subset of (one connected component of) the unframed moduli space 
$\widehat{\cm}(1,\Theta )$. $\Diamond$
\label{galler}
\end{theorem}

\noindent Theorem \ref{galler} can be regarded as a familiar ``collar 
theorem'' for the $1$-anti-instanton moduli space over the 
multi-Taub--NUT space. The patient reader may write down these 
solutions explicitly by inserting the explicit Green functions 
\cite[Equation (18)]{pag} of Page into (\ref{dekompozicio}) 
and then by the aid of the gauge (\ref{bazis}) construct the 
corresponding connection 1-form via (\ref{perturbacio}). However the 
quite complicated result is not informative. 

Rather we turn attention to a puzzling feature of the description 
provided by Theorem \ref{galler}. This is the limit 
$\lambda\rightarrow+\infty$. In this case (\ref{dekompozicio}) up to 
reparametrization remains meaningful and provides us with the positive 
minimal Green function $G(\cdot, y)$. It readily follows from the proofs of 
Lemmata \ref{mercetrafo} and \ref{aszimptotikus} 
that the resulting anti-self-dual connection 
$\nabla_{\tilde{A}_{y,+\infty}}$ gives rise to an unframed anti-instanton 
$[\nabla_{\tilde{A}_{y,+\infty}}]$ on 
$\Sigma^+$ obeying the weak holonomy condition with respect to 
$\nabla_{\Theta}$ and the rapid decay condition hence these solutions are 
expected to complete the moduli space. 

Indeed, for a fixed $y\in M_V$ taking the limit $\lambda\rightarrow 
+\infty$ (especially in (\ref{gorbulet})) the proof of Lemma 
\ref{mercetrafo} continues to hold 
hence the corresponding anti-self-dual connection is smooth 
everywhere on $\Sigma^+$. 

Additionally, in light of (\ref{alapfv}) the asymptotics of $G(\cdot, y)$ 
looks like $1/r$ therefore the proof of Lemma \ref{aszimptotikus} 
also can be repeated if we notice that this time $\vert a_y\vert$ 
and $\vert\dd_\Theta a_y\vert$ are $O(r^{-1})$ only. Hence 
$a_y\not\in L^2_{\frac{1}{2},1,\Theta}(M_V;\Lambda^1M_V\otimes{\rm 
End}\Sigma^+)$ but it follows from a direct calculation that the 
corresponding curvature still belongs to the weighted Sobolev 
space: $F_{\tilde{A}_y}\in L^2_{\frac{1}{2}}(M_V;\Lambda^2M_V\otimes{\rm 
End}\Sigma^+)$.
We conclude that $[\nabla_{\tilde{A}_{y,+\infty}}]$ provides us with an 
admissible unframed anti-instanton with unit energy on the 
positive spin bundle; more precisely it decays rapidly and one can 
also check that the 
weak holonomy condition with respect to the same trivial flat connection 
$\nabla_\Theta$ is fulfilled. If $y\in M_V$ is not a NUT then the 
corresponding anti-instanton is irreducible otherwise reducible to U$(1)$ 
and these Abelian solutions span the full $L^2$ cohomology of the 
multi-Taub--NUT geometry as it was stated in \cite[Theorem 5.1]{ete-hau3}. 

Therefore if $y\in M_V$ is not a NUT then these limiting 
irreducible solutions are expected to give the completion of the 
collar $M_V\times (0,+\infty )$ to the full unframed moduli space 
in the unconcentrated or ``centerless'' regime---or 
at least one connected component of the moduli space emerges this way. 
Since it is natural to include the reducible solutions corresponding to 
the NUTs into $\widehat{\cm}(1,\Theta )$ as well, for simplicity we 
shall continue to denote (one connected component of) this extended 
moduli space by $\widehat{\cm}(1,\Theta )$ from 
now on and note that it is not a manifold anymore: usual 
conical singularities appear around the reducible points as we know 
from the general theory but will also see shortly. 

We also conjecturially note that in our opinion not only a connected 
component but the whole moduli space of unit energy admissible 
anti-instantons with trivial holonomy arises this way.

\begin{theorem} Consider the multi-Taub--NUT space $(M_V, g_V)$ with $s>0$ 
NUTs $p_1,\dots,p_s\in M_V$, equipped with the natural orientation as 
above. Then for (one connected component of) the unframed moduli space 
of unit energy {\rm 
SU$(2)^+$} admissible anti-instantons decaying rapidly to the trivial 
flat connection on $\Sigma^+$ we find
\begin{equation}
\widehat{\cm}(1,\Theta )\cong (M_V\times (0,+\infty ])/\sim
\label{moduluster-eloallitas}
\end{equation}
where the equivalence relation $\sim$ means that $M_V\times\{ +\infty\}$ 
is pinched into $\R^3$ by collapsing the $S^1$-isometry orbits of $(M_V, 
g_V)$. 

Consequently there exists a singular fibration 
\begin{equation}
\Phi : \widehat{\cm}(1,\Theta )\longrightarrow\R^3
\label{moduluster-fibralas}
\end{equation}
with generic fibers homeomorphic to the open 
$2$-ball $B^2$ and as many as $s$ singular fibers 
homeomorphic to the semi-open $1$-ball $(0,+\infty ]$. Therefore (one 
connected component of) the 
moduli space is contractible and in particular is orientable. 

The images of the points $(p_i,+\infty )$ in $\widehat{\cm}(1,\Theta )$ 
with $i=1,2,\dots ,s$ represent reducible anti-instantons and 
$\widehat{\cm}(1,\Theta )$ around these points looks like a cone over 
$\C P^2$ or $\overline{\C P}^2$ (depending on the orientation).
\label{moduluster}
\end{theorem}

\noindent {\it Proof.} Athough we stated it already in Theorem
\ref{galler} first of all we remark that from the expression 
(\ref{gorbulet}) of the curvature one can immediately read off that both 
$y\in M_V$ and $\lambda\in (0,+\infty )$ are gauge-invariant parameters i.e., 
anti-instantons corresponding to different values of $(y,\lambda )$ cannot 
be gauge equivalent. Next, in Sect. \ref{five} we will prove that if 
$y\in M_V$ and ${\rm e}^{\ii\tau}y$ is its image by an element of 
the $S^1$ isometry group of $(M_V, g_V)$ then the corresponding two 
anti-instantons with $(y,\lambda =+\infty)$ and $({\rm e}^{\ii\tau}y, 
\lambda =+\infty )$ are actually gauge equivalent unit energy, rapidly 
decaying anti-instantons with trivial holonomy at infinity. 

In addition to these we know from \cite[Theorem 3.2]{ete-jar} that 
away from the reducible points the moduli space is a finite dimensional, 
possibly non-connected manifold without boundary. 

Consequently we are forced to look at our 
moduli space as follows. Fix a point $y\in M_V$ different from any NUT and 
consider the $S^1$ isometry orbit through $y$: it is parameterized by 
the cyclic coordinate $\tau\in [0,2\pi )$ of 
the local coordinate system on $(M_V, g_V)$ introduced in Sect. 
\ref{three}; as well as take the concentration parameter 
$\lambda\in (0,+\infty )$. Then by Theorem \ref{galler}, associated to 
any $y\in M_V$, we obtain a subset of the moduli space paremeterized by 
$(\tau ,\lambda )$. We regard this as a punctured 
open $B^2$ such that the points with $\lambda =0$ are supposed to 
represent the boundary while $\lambda =+\infty$ is the center of this 
$(\tau ,\lambda )$-ball. In this picture we interpret the existence of 
smooth solutions with $\lambda =+\infty$ as filling in the punctured 
$(\tau ,\lambda )$-ball over $y\in M_V$ with the anti-instanton 
$[\nabla_{\tilde{A}_{y,+\infty}}]$ 
corresponding to the limit $\lambda\rightarrow +\infty$. In other words, 
this anti-instanton represents the center of the unpunctured $(\tau 
,\lambda )$-ball that consequently will be homeomorphic to an open $B^2$ 
over any $y\in M_V$ different from the NUTs. 

Since the NUTs $p_1,\dots,p_s\in M_V$ are the only fixed points of the  
isometry group action on $(M_V, g_V)$ it follows that over these points
the moduli space is parameterized by $\lambda\in (0,+\infty ]$ only hence 
is a semi-open $B^1$. In this case the ``center'' is the endpoint 
$\lambda =+\infty$.

Therefore completing our collar $M_V\times (0,+\infty )$ 
with these limiting solutions, the moduli space takes the shape as in 
(\ref{moduluster-eloallitas}) where the equivalence relation $\sim$ 
means that $M_V\times\{+\infty\}$ is collapsed into $\R^3$ by shrinking all the 
isometry orbits of $(M_V, g_V)$ into single points respectively and admits a 
singular fibration (\ref{moduluster-fibralas}) whose generic fiber is an 
open $B^2$ and with $s$ pieces of singular fibers homeomorphic to a 
semi-open $B^1$. In particular $\widehat{\cm}(1,\Theta )$ is 
contractible 
hence orientable.

We already know that the only reducible solutions are 
the images $r_i\in\widehat{\cm}(1,\Theta )$ of the points $(p_i, 
+\infty )$ with $i=1,2,\dots,s$. We will study the structure of the moduli 
space around these points. Take a NUT $p_i\in M_V$ and consider its 
neighbourhood $U(r_i)\subset\widehat{\cm}(1,\Theta )$. We will construct 
a decomposition $\partial\overline{U}(r_i)=X_1\cup 
X_2$ of its boundary as follows. It is clear that if $y$ is not a NUT 
and the image of $(y,+\infty )$ in $\widehat{\cm}(1,\Theta )$ is 
$m$ then $U(m)\cong B^5$, an open 5-ball. Let us decompose its 
closure as $\overline{U}(m)\cong\overline{B}^3\times\overline{B}^2$ where 
$B^2$ is the aforementioned $(\tau ,\lambda )$-disk and $B^3$ is a ball 
around the image of $y$ in $\R^3$. Then $\partial\overline{U}(m)\cong 
(S^2\times\overline{B}^2)\cup (\overline{B}^3\times S^1)$ which is of 
course an $S^4$ glued together from $S^2\times\overline{B}^2$ and 
$\overline{B}^3\times S^1$ along their common boundary $S^2\times 
S^1$. Now let us move $y$ toward a particular NUT $p_i$; then one of the 
$\overline{B}^2$'s in the $S^2\times\overline{B}^2$ component of the 
above decomposition, namely that one which is moved exactly into the 
position of the NUT, gets shrink to a $\overline{B}^1$ in a 
manner such that meanwhile its radius (measured by $\lambda$) is kept 
constant, its circumference (measured by $\tau )$ gets 
vanish.\footnote{This way of collapsing $\overline{B}^2$ into 
$\overline{B}^1$ looks like closing an umbrella.} Put 
\[X_1\cong S^2\times\overline{B}^2/\sim_1\]
corresponding to this collapsed space. We obviously find for the singular 
cohomology that  
\[H^j(X_1;\Z )\cong\left\{\begin{array}{ll}
\Z & \mbox{if $j=0,2$}\\
 0 & \mbox{otherwise.}
\end{array}\right.\]
Since this movement also collapses an $S^1$, the boundary in the other
component $\overline{B}^3\times S^1$ of the previously collapsed 
$\overline{B}^2$, we also get a corresponding collapsed space
\[X_2\cong\overline{B}^3\times S^1/\sim_2 .\]
Since the generator of $\pi_1(\overline{B}^3\times S^1)\cong\Z$ is just 
killed out in $X_2$, we also obtain 
\[H^j(X_2;\Z )\cong\left\{\begin{array}{ll}
\Z & \mbox{if $j=0$}\\
 0 & \mbox{otherwise.}
\end{array}\right.\]
It is also clear that $X_1\cap X_2$ is nothing else than a $3$-sphere 
with exactly two antipodal points pinched consequently 
\[H^j(X_1\cap X_2;\Z )\cong\left\{\begin{array}{ll}
\Z & \mbox{if $j=0,3$}\\
 0 & \mbox{otherwise.}
\end{array}\right.\]

Recalling the Mayer--Vietoris sequence of the decomposition $\partial 
\overline{U}(r_i)=X_1\cup X_2$ which looks like
\[\dots\rightarrow H^j\left(\partial\overline{U}(r_i);\Z 
\right) \rightarrow H^j(X_1;\Z )\oplus H^j(X_2;\Z )\rightarrow 
H^j(X_1\cap X_2;\Z )\rightarrow H^{j+1}\left(\partial\overline{U}(r_i);\Z
\right)\rightarrow\dots\]
we end up with 
\[H^j\left(\partial\overline{U}(r_i);\Z\right)\cong
\left\{\begin{array}{ll}\Z & \mbox{if $j=0,2,4$}\\
                         0 & \mbox{otherwise.}
\end{array}\right.\]
Moreover we know that $\partial\overline{U}(r_i)$ is simply connected
and is smooth. However from all of these we obtain via Freedman's 
classification \cite{fre} that $\partial\overline{U}(r_i)\cong\C P^2$ or 
$\overline{\C P}^2$ depending on the orientation put onto 
$\widehat{\cm}(1,\Theta )$. 
Therefore the unframed moduli space about a reducible point looks like a 
cone over one of these spaces as stated. $\Diamond$ 
\vspace{0.1in}

\noindent We make a comment on the higher-energy 
moduli spaces including anti-instantons with trivial holonomy in 
infinity. It is obvious that instead of the functions
(\ref{dekompozicio}) we can use harmonic functions with a finite number 
of isolated singularities yielding

\begin{theorem}

The moduli spaces $\cm (k,\Theta )$ hence $\widehat{\cm}(k,\Theta 
)$ are non-empty for all $k\in\N$. $\Diamond$
\label{nagyenergia}
\end{theorem}

\noindent For sake of completeness without proof we exhibit here 
the degenerate $s=0$ case as well i.e., the flat 
$\R^3\times S^1$ whose associated instantons also referred to as {\it 
calorons}. This space is flat hence conformal rescaling works. Using 
coordinates $(x^1,x^2,x^3,\tau )=(x,\tau )$ on 
$\R^3\times S^1$ we have an obvious five paramerer family of positive 
harmonic functions (\ref{dekompozicio}) with Green functions 
in a rather explicit form \cite{pag}
\[G((x,\tau_x),(y,\tau_y))=\frac{1}{16\pi^2\vert 
x-y\vert}\frac{\tanh\vert x-y\vert}{1-\frac{\cos(\tau_x-\tau_y)}
{\cosh\vert x-y\vert}}=
\frac{1}{4\pi^2}\sum\limits_{k=-\infty}^{+\infty}\frac{1}{\vert 
x-y\vert^2+(\tau_x-\tau_y+2\pi k)^2}.\]
Making use of \cite[Theorem 2.3]{ete-jar} to describe the 
holonomy issue in this situation as well as \cite[Theorem 4.1]{ete-jar} 
which provides us with the dimension moreover taking into account that 
there are no reducible SU$(2)^\pm$ instantons on this space we state

\begin{theorem} Take the flat $\R^3\times S^1$ with an arbitrary fixed 
orientation. Then any non-trivial {\rm SU}$(2)^\pm$ instanton on 
$\Sigma^\pm$ is 
irreducible moreover the rapidly decaying ones obey the weak holonomy condition 
with respect to some flat connection $\nabla_\Gamma$ on the neck and 
their energies are always non-negative integers $k\in\N$. The 
corresponding framed moduli spaces, if not empty, have dimensions 
$\dim\cm (k,\Gamma )=8k$ while $\dim\widehat{\cm}(k,\Gamma )=8k-3$ for 
the unframed ones. 

Assume the asymptotic flat connection is the trivial one 
$\nabla_\Theta$ and let $\cm (1,\Theta )$ be the corresponding 
framed moduli space. Then for (one connected component of) the unframed 
space we find
\[\widehat{\cm}(1,\Theta )\cong (\R^3\times S^1\times (0,+\infty 
])/\sim\:\:\:\cong\R^3\times B^2\]
moreover the higher energy framed moduli spaces $\cm (k,\Theta 
)$ hence $\widehat{\cm}(k,\Theta )$ are also non-empty with $k\in\N$ 
arbitrary. $\Diamond$
\label{kaloron}
\end{theorem}


\section{Completion of the moduli space}
\label{five}

In this closing section our aim is to prove that conformal 
rescaling gives rise to at least one connected component of unit energy, 
rapidly decaying anti-instantons obeying the weak holonomy condition 
with respect to the trivial flat connection (however we note again that 
in our opinion this moduli space is connected). 

Our starting point is the construction of the twistor space of the 
multi-Taub--NUT geometry by the aid of Hitchin \cite{bes,hit1, hit3} and 
also a certain compactification of it. Then following Atiyah \cite{ati} 
this enables us to identify the harmonic functions used so far with 
elements of certain Ext groups on the compactified twistor space. 
These Ext groups also provide us with 
twisted holomorphic vector bundles over the original non-compact twistor 
space and the corresponding untwisted ones represent our 
anti-instantons in a new form in the spirit of the Atiyah--Ward 
correspondence. 

To get these twisted vector bundles more explicitly we 
investigate the real structure of the twistor space leading to the 
multi-Taub--NUT geometry and also construct the corresponding real 
twistor lines; our twisted vector bundles then, referring to 
Serre's method, can be constructed via sections vanishing along these real 
lines. The transition functions of these twisted bundles corresponding 
to pure Green functions, i.e., in the $\lambda\rightarrow 
+\infty$ limit can be constructed rather explicitly. We will 
find that the twisted bundles corresponding to $G(\cdot ,y)$ and 
$G(\cdot ,{\rm e}^{\ii\tau}y)$ are actually 
isomorphic hence the corresponding anti-instantons must be gauge 
equivalent verifying our picture on the moduli space presented in Theorem 
\ref{moduluster}. 

To begin with, we claim that there exist nice compactifications of the 
complex manifolds underlying a given a multi-Taub--NUT space.

\begin{lemma}
Consider any multi-Taub--NUT space $(M_V, g_V)$ with $s$ NUTs and let $X$ 
denote the smooth complex surface constructed in Lemma \ref{komplex} 
given by picking one generic complex structure in the hyper-K\"ahler family. 
Also let $\widetilde{X}$ be the smooth complex surface belonging to any 
exceptional complex structure, as constructed in Lemma \ref{komplex}. 

Then both $X$ and $\widetilde{X}$ admit compactifications 
$\overline{X}$ and $\overline{\widetilde{X}}$ respectively which are 
smooth compact rational surfaces. These spaces arise by attaching $s+3$ 
lines to the finite parts in a suitable way.
\label{Xkompaktifikalas} 
\end{lemma}

\noindent{\it Proof.} Fix a generic configuration of $s$ points 
$q_1,\dots,q_s\in\R^3$ 
and consider the corresponding multi-Taub--NUT space $(M_V, g_V)$. Pick 
a direction $e_1$ in $\R^3$ not parallel with any 
straight line segment $\ell_{ij}$ connecting $q_i$ and $q_j$ and put the 
corresponding (generic) integrable complex structure $J_{e_1}$ onto 
$M_V$ from the hyper-K\"ahler family. We know from Lemma \ref{komplex} 
that $(M_V, J_{e_1})$ is biholomorphic to the algebraic surface 
$X\subset\C^3$ given by (\ref{model}) and there are $s(s-1)$ exceptional 
directions to choose for $e_1$ so that the corresponding complex manifolds 
are {\it not} biholomorphic to $X$ but rather to its blowup denoted by 
$\widetilde{X}$ in agreement with the notations of Lemma \ref{komplex}.

First consider the case of $X$. Introducing homogeneous coordinates 
$([x:u],[y:v],[z:w])\in\C P^1_x\times\C P^1_y\times\C P^1_z$ equation 
(\ref{model}) can be made homogeneous of degrees $(1,1,s)$ 
respectively as follows:
\begin{equation}
xyw^s-uv\prod_{i=1}^s(z-wp_i)=0.
\label{projektiv}
\end{equation}
Let us denote the resulting compact complex surface in $\C P^1_x\times\C 
P^1_y\times\C P^1_z$ by $\overline{X}^*$. We claim that $\overline{X}^*$ 
is smooth except one point. Indeed, let $0_x:=[0:1]$ and $\infty_x:=[1:0]$ 
denote the origin and the infinity of $\C P^1_x$ respectively and in the same 
fashion introduce the notation $0_y ,\infty_y\in\C P^1_y$ and 
$0_z,\infty_z\in\C P^1_z$. Then in the compactification 
$\C^3\subset\C P^1_x\times\C P^1_y\times\C P^1_z$ infinity is 
represented by the bouquet
\[\infty :=(\{\infty_x\}\times\C P^1_y\times\C P^1_z)\vee (\C 
P^1_x\times\{ \infty_y\}\times\C P^1_y)\vee (\C P^1_x\times\C 
P^1_y\times\{\infty_z\} )\]
and one obtains that $\overline{X}^*\cap\infty$ consists of four lines, 
that is $\overline{X}^*=X\sqcup (\ell_1\cup \ell_2\cup \ell_3\cup\ell_4)$ with 
\[\ell_1:=\{\infty_x\}\times\{ 0_y\}\times\C 
P^1_z,\:\:\:\:\:\ell_2:=\{ \infty_x\}\times\C P^1_y\times\{ 
\infty_z\},\]
\[\ell_3:=\C P^1_x\times\{\infty_y\}\times\{\infty_z\},\:\:\:\:\:
\ell_4:=\{ 0_x\}\times\{\infty_y\}\times\C P^1_z.\]
These lines are not disjoint, they intersect each other as follows: 
\[\ell_1\cap\ell_2=(\infty_x,0_y,\infty_z),
\:\:\:\:\:\ell_2\cap\ell_3=(\infty_x,\infty_y,\infty_z),\:\:\:\:\:\ell_3\cap 
\ell_4=(0_x,\infty_y,\infty_z).\]
By calculating the gradient of the left hand side in (\ref{projektiv}) 
we see that the only point where $\overline{X}^*$ is singular is the 
intersection $\ell_2\cap\ell_3$ where $\overline{X}^*$ locally looks like
\begin{equation}
uv-\prod\limits_{i=1}^sw(1-wp_i)^{-1}=0
\label{vegtelen-szingularitas}
\end{equation}
i.e., it possesses a classical $A_{s-1}$ quotient singularity.

Carrying out the same procedure for the exceptional singular surface 
$X^*\subset\C^3$ with the notations of Lemma \ref{komplex}, we obtain a 
similar compactification $\overline{X^*}^*$ in 
$\C P^1_x\times\C P^1_y\times\C P^1_z$. The only difference is that in 
addition to the $A_{s-1}$ singularity at 
infinity as above, $\overline{X^*}^*$  possess more rational 
double points in the finite part in each point where $[p_i:1]=[p_j:1]$.

We proceed further and desingularize these singular compactifications. In 
light of Lemma \ref{komplex} we simply do this by replacing them with 
their non-singular blowups. Again consider first the generic case 
$\overline{X}^*$ and let us denote the
resulting smooth space by $\overline{X}$. Clearly, in this case we have to 
just remove the $A_{s-1}$-singularity at infinity from 
$\overline{X}^*$ therefore we have to attach $4+s-1=s+3$ lines in an 
approriate way:
\begin{equation}
\overline{X}=X\sqcup (\ell_1\cup\ell_2\cup\ell '_1\cup\dots
\cup\ell '_{s-1}\cup\ell_3\cup\ell_4)
\label{felfujas}
\end{equation}
to obtain $\overline{X}$. Denoting by $\overline{\widetilde{X}}$ 
the desingularization of $\overline{X^*}^*$, it stems 
exactly the same way from the smooth $\widetilde{X}$, the 
desingularization of $X^*\subset\C^3$. 

Finally we prove that both of these smooth compactifications are rational 
surfaces. Indeed, consider first general case when $\overline{X}^*$ is 
smooth except at infinity. Then $\overline{X}^*$ is given by equation
(\ref{projektiv}) for distinct values of $p_j$'s. Consider the map
\[\overline{X}^*\longrightarrow\C P^1_x\times\C P^1_z,\:\:\:\:\:
    ([x:u],[ y:v],[z:w])\longmapsto ([x:u],[z:w]).\]
Over a given $([x:u],[z:w])\in \C P^1_x\times\C P^1_z$ there exists a  
unique $[y:v]\in  \C P^1_y$ such that $([x:u],[y:v],[z:w])\in
\overline{X}^*$, except in the case $xw^s=u\prod (z-wp_j)=0$ when
any $[y:v]$ is a solution. It readily follows that these exceptional
points are precisely
\begin{equation}
\{ (0_x,[p_1:1]),\ldots ,(0_x,[p_s:1]),(\infty_x ,\infty_z )\}  
\in\C P^1_x\times\C P^1_z.
\label{exceptionalpoints}
\end{equation}
Furthermore, all the points of $\overline{X}^*$ are non-singular except
the isolated hence codimension $2$ surface singularity of type $A_{s-1}$
at $(\infty_x ,\infty_y ,\infty_z )\in\C P^1_x\times\C P^1_y\times\C 
P^1_z$ which is mapped onto $(\infty_x,\infty_z)$ in $\C
P^1_x\times\C P^1_z$. It follows that $\overline{X}^*$ is
isomorphic to the blowup of $\C P^1_x\times\C P^1_z$ in the points
(\ref{exceptionalpoints}) outside of a codimension $2$ subset.  
The same argument works for the special case of $\overline{X^*}^*$ 
possessing a further rational double point, if we
add this singular point to the codimension $2$ subset above.
Therefore both $\overline{X}^*$ and $\overline{X^*}^*$ are rational. Since 
both $\overline{X}$ and $\overline{\widetilde{X}}$
arise by blowing further up $\overline{X}^*$ and $\overline{X^*}^*$ 
respectively in finitely many
points, we conclude that the non-singular compact spaces $\overline{X}$ 
and $\overline{\widetilde{X}}$ are also rational. 
$\Diamond$
\vspace{0.1in}

\begin{remark}\rm It is worth comparing the complex compactification
(\ref{felfujas}) with the simple smooth real compactification
of $M_V$ in \cite{ete2,ete-jar} motivated by $L^2$-cohomology theory of 
ALF spaces \cite{hau-hun-maz}.
\end{remark}

\noindent We proceed further and smoothly compactify the 
twistor space of a multi-Taub--NUT space in two steps.

\begin{lemma}
The {\em first approximation of the twistor space} (what we denote here 
by) $Z$ of any multi-Taub--NUT space $(M_V, g_V)$ with $s$ NUTs admits a 
smooth complex compactification (what we denote here by) $\ZZ$. This 
compactification arises by adding finitely many various Hirzebruch surfaces 
to this $Z$ in a suitable way.
\label{Zkompaktifikalas1}
\end{lemma}

\noindent{\it Proof.} First we construct a sort of 
approximation of the twistor space of $(M_V, g_V)$ 
and then modify it into the true twistor space in the next lemma.

Being $(M_V, g_V)$ simply connected with 
vanishing Ricci curvature, its true twistor space admits a holomorphic 
fibration over $\C P^1$ consequently we would like to regard equation 
(\ref{projektiv}) as a set of equations parameterized by a projective 
line. Therefore we proceed as follows (cf. \cite[pp. 
393-395]{bes} and \cite{hit1}). Consider the holomorphic fiber bundle 
\begin{equation}
\overline{\pi} :P(H^k\oplus H^0)\oplus P(H^l\oplus H^0)\oplus 
P(H^m\oplus H^0)\longrightarrow\C P^1
\label{proj.fibralas}
\end{equation}
whose fibers are $F_x\times F_y\times F_y$ with $F_i\cong\C P^1_i$ 
($i=x,y,z$) as before. Take the subbundle $\overline{\pi}_x:P(H^k\oplus 
H^0)\rightarrow \C P^1$ with $F_x\cong\C P^1_x$. Referring to the isomorphism 
\[P(H^0\oplus H^{-k})\stackrel{\otimes H^k}{\longrightarrow} 
P(H^k\oplus H^0)\]
take the canonical section $\infty_x:=(1,0)$ of $P(H^0\oplus H^{-k})$ and 
regard its image as a divisor which we also denote by $\infty_x$; 
one then has an associated line bundle $L_x:=[\infty_x]$ and 
corresponding sheaf $\co_{{\rm rel},x}(1)$ on $P(H^0\oplus H^{-k})$. 
This line bundle has a canonical section $u\in \co_{{\rm rel},x}(1)$ 
satisfying $(u)=\infty_x$. We also have a canonical section 
$0_x:=(0,1)$ of $P(H^k\oplus H^0)$ which gives 
rise to the bundle $L_x\otimes\overline{\pi}^*_xH^k$ and a canonical 
section $x$ which is an element of $\co_{{\rm rel},x}(1)\otimes\co_{P(H^k\oplus 
H^0)}(k)$ where $\co_{P(H^k\oplus H^0)}(k)$ is the associated sheaf of 
$\overline{\pi}^*_xH^k$. Take the fibration
\[\overline{\pi}_{y,z}:P(H^k\oplus H^0)\oplus P(H^l\oplus
H^0)\oplus P(H^m\oplus H^0)\longrightarrow P(H^k\oplus H^0)\]
whose fibers are $F_y\times F_z$. Pulling back $x,u$ we therefore obtain 
elements 
\[\overline{\pi}^*_{y,z}x\in\overline{\pi}^*_{y,z}\co_{{\rm
rel},x}(1)\otimes\overline{\pi}^*_{y,z}\co_{P(H^k\oplus 
H^0)}(k),\:\:\:\:\:\overline{\pi}^*_{y,z}u\in\overline{\pi}^*_{y,z}\co_{{\rm 
rel},x}(1).\]
In the same fashion we construct 
\[\overline{\pi}^*_{x,z}y\in\overline{\pi}^*_{x,z}\co_{{\rm
rel},y}(1)\otimes\overline{\pi}^*_{x,z}\co_{P(H^l\oplus 
H^0)}(l),\:\:\:\:\:\overline{\pi}^*_{x,z}v\in\overline{\pi}^*_{x,z}\co_{{\rm 
rel},y}(1)\] 
and finally 
\[\overline{\pi}^*_{x,y}z\in\overline{\pi}^*_{x,y}\co_{{\rm
rel},z}(1)\otimes\overline{\pi}^*_{x,y}\co_{P(H^m\oplus 
H^0)}(m),\:\:\:\:\:\overline{\pi}^*_{x,y}w\in\overline{\pi}^*_{x,y}\co_{{\rm 
rel},z}(1).\] 
For simplicity we shall denote them as $x,u,y,v,z,w$ respectively and 
observe that they provide us with coordinates on the total space of the 
fibration (\ref{proj.fibralas}).

The situation gets simplified somewhat by virtue of the following 
observation. It is not difficult to check that the 
bundle $\overline{\pi}^*_{y,z}(\overline{\pi}^*_xH^k)$ 
on the total space of (\ref{proj.fibralas}) depends only on $k$ and not 
on the way we pulled it back from $\C P^1$. Consequently
\[\overline{\pi}^*_{y,z}(\overline{\pi}^*_xH^k)\cong\overline{\pi}
^*_{x,y}(\overline{\pi}^*_zH^k)
\cong\overline{\pi}^*_{x,z}(\overline{\pi}^*_yH^k)\cong\overline{\pi}^*H^k\]
and the same is true for the other bundles $H^l$ and $H^m$. Write  
\[W:=P(H^k\oplus H^0)\oplus P(H^l\oplus H^0)\oplus P(H^m\oplus H^0)\]
for the total space of (\ref{proj.fibralas}) and put 
\[L^{a_0,a_1,a_2,a_3}:=\overline{\pi}^*H^{a_0}\otimes\overline{\pi}
^*_{y,z}L^{a_1}_x\otimes
\overline{\pi}^*_{x,z}L^{a_2}_y \otimes\overline{\pi}^*_{x,y}L^{a_3}_z\]
with $a_i\in\Z$; we will also use the convenient notation for the 
corresponding sheaf:
\[\co_W(a_0,a_1,a_2,a_3):=\co_W(a_0)\otimes\overline{\pi}^*_{y,z}\co_{{\rm 
rel},x}(a_1)\otimes\overline{\pi}^*_{x,z}\co_{{\rm 
rel},y}(a_2)\otimes\overline{\pi}^*_{x,y}\co_{{\rm rel},z}(a_3)\] 
(therefore $\co_W(a_0)=\co_W(a_0,0,0,0)$). Then $x\in\co_W(k,1,0,0)$ 
and $u\in \co_W(0,1,0,0)$, etc. and if we regard $p_i\in\co_W(m,0,0,0)$ 
we obtain 
\[xyw^s\in\co_W(k+l, 1,1,s)\:\:\:\:\:{\rm and}\:\:\:\:\:
uv\prod\limits_{i=1}^s(z-wp_i)\in\co_W(sm,1,1,s)\]
consequently setting $k+l=sm$ we end up with a well-defined map
$\overline{P}: W\longrightarrow L^{sm,1,1,s}$ where formally 
$\overline{P}$ is given by the 
left hand side of (\ref{projektiv}) and a singular model for the first 
approximation of the compactified twistor space $\ZZ^*\subset W$ is provided 
by the hypersurface $\overline{P}=0$ hence restricting (\ref{proj.fibralas}) 
it admits a holomorphic fibration $\overline{\pi}: \ZZ^*\longrightarrow\C P^1$.

To simplify the notation of Lemma \ref{Xkompaktifikalas} from now on 
we shall denote the fibers $\overline{\pi}^{-1}([a:b])$ of 
this fibration as $\overline{X}^*_{[a:b]}$ for all  
$[a:b]\in\C P^1$. Note that $\ZZ^*$ has an obvious singularity along the 
line $\ell$ which is the section passing through the singular points 
$(\infty_x,\infty_y,\infty_z)_{[a:b]}\in\overline{X}^*_{[a:b]}$ for an 
$[a:b]\in\C P^1$. However $\ZZ^*$ has more singularities, namely the 
points $(0, 0, p_i([a:b]))$ where exactly two 
roots coincide: $p_i([a:b])=p_j([a:b])$. (We may assume that the 
configuration is generic hence no more than two roots coincide at the 
same time.) Hence these points represent rational 
double points in $\ZZ^*$. Therefore the number of such bad points is 
$m\cdot\binom{s}{2}=m\cdot\frac{1}{2}s(s-1)$ since $p_i\in 
\co_W(m,0,0,0)$. 
The geometric origin of these bad points is clear: they just correspond to 
the special complex structures analyzed in Lemma \ref{komplex} 
under which $(M_V, J_{e_1})$ is {\it not} 
biholomorphic to (\ref{model}). Since the number of these complex structures is 
$s(s-1)$ as we have seen, we have to set $m=2$ yielding $k+l=2s$ 
and since $x$ and $y$ play symmetric roles, we take $k=l=s$ in 
(\ref{proj.fibralas}). 

We desingularize $\ZZ^*$ to obtain a smooth compact complex 
$3$-manifold $\ZZ$ whose finite part $Z\subset\ZZ$ provides us with a 
non-singular model for the approximated twistor space of $(M_V, g_V)$. In 
light of Lemmata \ref{komplex} and \ref{Xkompaktifikalas} we simply do this 
by replacing the fibers $\overline{X}^*_{[a:b]}$ of $\overline{\pi}: 
\ZZ^*\longrightarrow\C P^1$ with their 
non-singular blowups for all $[a:b]\in\C P^1$. Let us denote the 
resulting smooth fiber by $\overline{X}_{[a:b]}$. Since these fibers 
are diffeomorphic to each other via Lemma \ref{komplex} we end up with a 
smooth holomorphic fibration; we denote it as
$\overline{\pi} :\ZZ\longrightarrow\C P^1$.

Finally, referring to (\ref{felfujas}) we conclude that the 
compactified smooth $\ZZ$ arises by adding $s+3$ Hirzebruch surfaces to 
the smooth twistor space $Z$ of $(M_V, g_V)$ as follows:
\[\ZZ =Z\sqcup (S_1\cup S_2\cup S'_1\cup\dots\cup S'_{s-1}\cup S_3\cup S_4).\]
Two of them, $S_1$ and $S_4$ stemming from the lines $\ell_1$ and
$\ell_4$ in $\overline{X}^*$, are biholomorphic to $P(H^2\oplus H^0)$ 
while $S_2$ and $S_3$ given by $\ell_2$ and $\ell_3$
are biholomorphic to $P(H^s\oplus H^0)$.
To determine the type of the exceptional ones, notice that the local
resolution of the singularity (\ref{vegtelen-szingularitas}) is given
by defining the new variables $u_s:=w^{-s}u,v_s:=v,w_s:=w$,
whence the equation of the proper transform becomes
\[   u_sv_s - \prod\limits_{i=1}^s(1-w_sp_i)^{-1}=0.\]
This means that we add a family of projective lines parametrized by
$\C P^1$ with sections $w_s$. As $w_s$ is untwisted in the base 
direction, this implies that the Hirzebruch surfaces 
$S'_1,\dots,S'_{s-1}$ we add
through this process are trivial $\C P^1 \times \C P^1$'s. $\Diamond$
\vspace{0.1in}

\begin{remark}\rm The reader may recognize that this ``first 
approximation'' of the twistor space of the multi-Taub--NUT space is 
nothing else than the twistor space of the corresponding ALE 
Gibbons--Hawking space \cite{hit1}.
\end{remark}

\begin{lemma} The {\em true twistor space} $Z$ of any multi-Taub--NUT space 
$(M_V, g_V)$ with $s$ NUTs also admits a smooth complex 
compactification $\overline{Z}$. It also arises by adding $s+3$ Hirzebruch 
surfaces to $Z$ in a suitable way.
\label{Zkompaktifikalas2}
\end{lemma}

\noindent{\it Proof.} The idea is to twist (cf. \cite[pp.
393-395]{bes} and \cite{hit3}) the natural algebraic model 
constructed in Lemma \ref{Zkompaktifikalas1} with a section of a certain 
line bundle $L^c$ which lives on the total space of the line 
bundle $\pi_z:H^2\rightarrow\C P^1$ where this base is the same 
as that of (\ref{proj.fibralas}).

For the definition of $L^c$, let this $\C P^1$ be covered by the 
affine open sets $U_a$ and $U_b$ with coordinates $[a:b]$ satisfying 
$a\not= 0$ and $b\not=  0$ respectively and write $U'_a=\pi_z^{-1}(U_a)$ and   
$U'_b=\pi_z^{-1}(U_b)$ for the open subsets of $H^2$ covering the affine 
subsets of $\C P^1$ moreover denote by $z$ any point in the total space of 
$H^2$. Therefore we obtain coordinate systems $(\left[ 1:\frac{b}{a} 
\right],z)$ and $(\left[ \frac{a}{b}:1\right] ,z)$ on $U'_a$ and $U'_b$
respectively. Define the holomorphic line bundle $L^c$ with a complex 
parameter $c$ by the transition function
\[g:U'_a\cap U'_b\longrightarrow\C^*,\:\:\:\:\:g(z,[a:b]):={\rm e}^{-c
z\frac{b}{a}} .\]
It is clear that $L^0$ is just the trivial bundle over $H^2$ and 
$L^{-c}$ is canonically isomorphic to $(L^c)^*$.

Next we construct a section of $L^c$. For this we exploit a 
nice geometric description of $H^2\cong T\C P^1$ by identifying it with 
the space of oriented affine lines in $\R^3$ (with a fixed origin) as 
follows. If $\ell\subset\R^3$
is an oriented affine line and $\ell (t)\in\ell$ is its point, $t\in\R$,
then clearly exist unique vectors ${\bf u}, {\bf v}$ with $\vert{\bf
u}\vert =1$ and ${\bf u}\cdot{\bf v}=0$ such that $\ell (t)={\bf
v}+t{\bf u}$. Regarding ${\bf u}$ as a point on $S^2$ and ${\bf v}$ a
vector in its tangent space we obtain the isomorphism. There is a natural
real structure on the space of oriented real lines
namely $({\bf u}, {\bf v})\mapsto (-{\bf u}, {\bf v})$ which gives rise
to the antipodal map on $\C P^1$ inside $H^2$ as a zero section. 
Consequently we also obtain a map from the space of
holomorphic sections into that of anti-holomorphic sections 
\begin{equation}
r: H^0(\C P^1;\co (2))\longrightarrow \overline{H^0}(\C P^1;\co (2)). 
\label{valos}
\end{equation}
A section is called {\it real} if it is invariant under this map. 

If $z\in H^2$ then there is a unique oriented affine line 
$\ell_z\subset\R^3$ corresponding to it 
with a distinguished point $\ell_z(0)\in\R^3$ (depending on the choice of 
origin in $\R^3$). This point can be characterized by the set of all lines 
in $\R^3$ intersecting it which provides us with a real section in $H^0(\C
P^1;\co (2))\cong\C^3$ hence is of the form $z([a:b])=\alpha a^2+2\beta ab
-\overline{\alpha}b^2$ with $\alpha\in\C$ and $\beta\in\R$.
Therefore we have two ways of looking at the same object: we can
regard $z$ either as a point on the total space $z\in H^2$ or as a real
section $\alpha a^2+2\beta ab
-\overline{\alpha}b^2\in H^0(\C P^1;\co (2))$.
Of course, quite tautologically speaking, the section $z$ passes
through the point $z$ corresponding to it.

Regarding $z$ as a real section $z([a:b])=\alpha a^2+2\beta ab
-\overline{\alpha}b^2$ now, this
homogeneous polynomial is canonically indentified with
the single-valued polynomial in the affine coordinate $\frac{a}{b}$ over $U_b$:
\[\left. z\right\vert_{U_b}\left(\frac{a}{b}\right) =
\alpha\left(\frac{a}{b}\right)^2+2\beta
   \left(\frac{a}{b}\right) -\overline{\alpha}. \]
Now we define a nowhere vanishing section $\eta^c$ of $L^c$ by
\[\eta^c (z,[a:b]):=\left\{\begin{array}{ll}
                                  {\rm e}^{c
\left(\alpha\left(\frac{a}{b}\right)+\beta\right)} & \mbox{on
$U'_b$;}\\
g(z,[a:b]){\rm e}^{c
\left(\alpha\left(\frac{a}{b}\right)+\beta\right)} &
\mbox{on $U'_a$.}
                              \end{array}\right.\]
It is well-defined since it looks like
$\exp c\left(\alpha\left(\frac{a}{b}\right)+\beta\right)$ on
$U'_b$ i.e., where $b\not= 0$ and looks like
$\exp c\left(\overline{\alpha}\left(\frac{b}{a}\right)-\beta\right)$
on $U'_a$ i.e., where $a\not= 0$. 

Next we claim that the pair $(L^c,\eta^c)$ can be extended uniquely to a 
pair $(\overline{L}^c, \overline{\eta}^c)$ where $\overline{L}^c$ is a 
line bundle over $\overline{\pi}_z: P(H^2\oplus 
H^0)\rightarrow\C P^1$ which restricts to $L^c$ over the finite part as 
before moreover 
\[\overline{\eta}^c\in H^0\left( P(H^2\oplus H^0)\:;\: 
\ce_{P(H^2\oplus H^0)}(\overline{L}^c)\right)\]
is a section of it which also gives back $\eta^c$ on $L^c$.

In order to construct these extensions we proceed as 
follows. First, notice that there is a covering of $P(H^2\oplus H^0)$ by 
four affine charts: the open sets $U'_b$ and $U'_a$ of the total space of 
$H^2$ defined above and additional charts $V'_b$ and $V'_a$ covering 
the section at infinity $(w=0)$. Explicitly, let $(z,a)$ and $(z',b)$ be 
standard coordinates on $U'_b$ and $U'_a$ respectively with $b=a^{-1}$ and 
$z'=za^{-2}$, then coordinates on $V'_b$ and $V'_a$ can be chosen to be 
$(w,a)$ and $(w,b)$ with $w=z^{-1}$. Since we want to 
extend the bundle $L^c$, the transition function of
$\overline{L}^c$ from the chart $U'_b$ to $U'_a$ must be that of $L^c$,
i.e. equal to $g=\exp(-cz/a)$. On the other hand, the requirement
that $\eta^c$ should extend to a section of $\overline{L}^c$
fixes the transition function between the charts $U'_b$ and $V'_b$ to be
$\exp(-c \alpha a - c \beta)$. Similarly, the transition between $U'_a$
and $V'_a$ is constrained to be $\exp(c \overline{\alpha} b + c \beta)$. 
From the cocycle-relation we then deduce that the restriction to the 
total space of $H^2$ (i.e. the open part of $P(H^2\oplus H^0)$ away from 
the section at infinity) of the transition between $V'_b$ and $V'_a$ 
must be identically $1$. In particular, this transition function extends 
to infinity and the set of transition functions just described yields the 
line bundle $\overline{L}^c$ with a section $\overline{\eta}^c$ we were 
looking for.

Then referring to the notations of Lemma \ref{Zkompaktifikalas1} (and 
knowing already that $k=l=s$ and $m=2$ as before) consider the fibration
\[\overline{\pi}_{x,y}:W\longrightarrow P(H^2\oplus H^0)\]
and use it to pull back both $\overline{L}^c$ and $\overline{\eta}^c$ over 
$W$. As an extension of the construction in Lemma \ref{Zkompaktifikalas1} 
set
\[L^{a_0,a_1,a_2,a_3,c}:=L^{a_0,a_1,a_2,a_3}\otimes\overline{\pi}^*_{x,y}
\overline{L}^c\]
with $a_i\in\Z$ and $c\in\C$; we will also use the convenient notation for 
the corresponding sheaf:
\[\ce_W(a_0,a_1,a_2,a_3,c):=\co_W(a_0,a_1,a_2,a_3)
\otimes\ce_W\left(\overline{\pi}^*_{x,y}\overline{L}^c\right) .\]
If---again in analogy with Lemma \ref{Zkompaktifikalas1}---we denote 
$\overline{\pi}^*_{x,y}\overline{\eta}^{\pm c}\in 
\ce_W\left(\overline{\pi}^*_{x,y}\overline{L}^{\pm c}\right)$ simply as 
$\overline{\eta}^{\pm c}$ then 
we have $x\cdot\overline{\eta}^c\in\ce_W(k,1,0,0,c)$ 
and $y\cdot\overline{\eta}^{-c}\in\ce_W(k,1,0,0,-c)$ etc., and can 
regard them as twisted coordinates on $W$. Making use of these new twisted 
coordinates in Lemma \ref{Zkompaktifikalas1} we 
define the {\it true compactified singular twistor space} $\ZZ^*$ of the 
multi-Taub--NUT geometry by the equation $\overline{P}_c=0$ where 
$\overline{P}_c: W\rightarrow L^{2s,1,1,s,0}$ looks like
\[\overline{P}_c(x,u,y,v,z,w):=
(x\cdot\overline{\eta}^c)(y\cdot\overline{\eta}^{-c})w^s-
(u\cdot\overline{\eta}^c)(v\cdot\overline{\eta}^{-c})
\prod\limits_{i=1}^s(z-wp_i)\]
which simply gives back (\ref{projektiv}). Therefore we come up 
with a singular holomorphic fibration
\begin{equation}
\overline{\pi}: \ZZ^*\longrightarrow\C P^1.
\label{szing.tvisztor.fibralas}
\end{equation}
Carrying out desingularization of $\ZZ^*$ as before we obtain finally that 
there is a smooth complex fibration
\begin{equation}
\overline{\pi} :\ZZ\longrightarrow\C P^1
\label{tvisztor.fibralas}
\end{equation}
whose fibers $\overline{X}_{[a:b]}$ are biholomorphic to
the rational surfaces constructed in Lemma \ref{Xkompaktifikalas}. 
This is the {\it true compact non-singular twistor space} of 
$(M_V, g_V)$. 

For clarity we remark that from now on, the notations $\ZZ^*$, $\ZZ$, 
$Z^*$ and $Z$ will denote the corresponding {\it true} 
twistor spaces of the multi-Taub--NUT geometry (and we will call them 
simply as the {\it twistor space} of the corresponding kind). $\Diamond$
\vspace{0.1in}

\noindent Now we are in a position to read off certain cohomology groups 
of the compactified twistor space.

\begin{lemma}
The middle sheaf cohomology groups of the compactified smooth twistor 
space $\ZZ$ constructed in Lemma \ref{Zkompaktifikalas2} satisfy
\[H^1(\ZZ ;\co_{\ZZ}(-2))\cong\C\:\:\:\:\:{\it and}\:\:\:\:\: H^2(\ZZ
;\co_{\ZZ}(-2))\cong 0.\]
\label{kohomologia}
\end{lemma}

\noindent{\it Proof.} First note that in light of Lemma 
\ref{Xkompaktifikalas} all the fibers $\overline{X}_{[a:b]}$ with 
$[a:b]\in\C P^1$ of (\ref{tvisztor.fibralas}) are rational.

Next, notice that $\co_{\ZZ}(-2)=\overline{\pi}^*\co (-2)$ is trivial 
along the fibers of (\ref{tvisztor.fibralas}). Rationality 
implies in particular that for all $[a:b]$ and $1\leq q\leq 2$ one has 
\[H^q\left(\overline{X}_{[a:b]};(\overline{\pi}\vert_{\overline{X}_{[a:b]}})
^*\co (-2)\right)\cong H^q\left(\overline{X}_{[a:b]};
\co_{\overline{X}_{[a:b]}}\right) =0,\]
so $R^q\overline{\pi}_*\overline{\pi}^*\co (-2)$ is the zero sheaf. On 
the other hand $R^0\overline{\pi}_*\overline{\pi}^*\co (-2) \cong\co 
(-2)$ by the projection formula consequently 
\[H^0\left(\C P^1;R^0\overline{\pi}_*\overline{\pi}^*\co(-2)\right) 
=0\:\:\:\:\:{\rm and}\:\:\:\:\:H^1\left(\C 
P^1;R^0\overline{\pi}_*\overline{\pi}^*\co(-2)\right) =\C .\] 
The second level of the Leray spectral sequence associated to 
(\ref{tvisztor.fibralas}) and the sheaf $\overline{\pi}^*\co(-2)$ is 
given by  $E^{p,q}_2=H^p\left(\C P^1; 
R^q\overline{\pi}_*\overline{\pi}^*\co(-2)\right)$. Taking 
into account our calculations so far, the only non-trivial part is 
a one-dimensional space in bidegree $(1,0)$: 
\[{}^q\underset{p}{\underline{\left\vert
         \begin{array}{lll}
               0 & 0 \\
               0 & 0 \\
               0 & \C   
         \end{array}
\right.}}\] 
implying that the spectral sequence collapses at this step and the statement 
follows. $\Diamond$
\vspace{0.1in}

\noindent Keeping these results in mind now we recall the notion of the 
``global Ext'' groups $\ext^k(V; \crr ,\cs )$ over a complex 
manifold $V$ with two coherent sheaves $\crr$ and $\cs$ on it. 
The background material can be found in \cite{gri-har}. These objects are 
appropriate generalizations of the sheaf cohomology groups over complex 
manifolds. For instance if $\crr =\co_V$ then we know that 
\begin{equation}
\ext^k(V; {\co_V},\cs )\cong H^k(V; \cs ).
\label{ext-H-atteres}
\end{equation}
If $V$ is compact $n$ dimensional then there is 
a generalization of the Serre duality theorem for coherent sheaves, 
called the {\it Grothendieck duality theorem}:
\begin{equation}
H^k(V;\cs )\cong (\ext^{n-k}(V;\cs ,\ck_V))^*
\label{grothendieck}
\end{equation}
and in the same fashion if $\cs =\co (E)$ is locally free then the 
{\it Serre duality theorem} implies that
\begin{equation}
(\ext^{n-k}(V;\co (E),\ck_V))^*\cong H^{n-k}(V;\co (E^*)\otimes\ck_V).
\label{serre}
\end{equation}
We will be particularly interested in $\ext^1(V; \crr ,\cs )$ 
which classifies extensions of coherent sheaves over $V$:
\[0\longrightarrow\cs\longrightarrow\ce\longrightarrow\crr\longrightarrow 
0.\]
Following Atiyah \cite{ati} we obtain:

\begin{lemma} Let $\cj_Y$ denote the ideal sheaf of 
holomorphic functions on $\ZZ$ vanishing on a fixed real line $Y\subset 
Z\subset\ZZ$. Then 
\[\ext^1(\ZZ;\cj_Y,\co_{\ZZ}(-2))\cong\C^2.\]
\label{ext-lemma}
\end{lemma}
\noindent{\it Proof.} Note that $\co_Y(-2)\cong\ck_Y$, the canonical 
sheaf of $Y$. Let $j:Y\to \ZZ$ be the inclusion and denote by $j_*\ck_Y$ 
the extension of $\ck_Y$ by zero over $\ZZ$. Hence $j_*\ck_Y$ is a 
coherent sheaf on $\ZZ$. Applying (\ref{grothendieck}) on $\ZZ$ in the 
last step we obtain
\[H^k(Y;\ck_Y)\cong H^k(\ZZ ;j_*\ck_Y)\cong 
(\ext^{3-k}(\ZZ ;j_*\ck_Y,\ck_{\ZZ}))^*\]
and additionally we know that 
\[\ext^{3-k}(\ZZ; j_*\ck_Y,\ck_{\ZZ})\cong 
\ext^{3-k}(\ZZ; j_*\co_Y \otimes \overline{\pi}^*\ck_Y,\ck_{\ZZ})
\cong 
\ext^{3-k}(\ZZ;j_*\co_Y,\overline{\pi}^*\ck^*_Y\otimes\ck_{\ZZ}).\]
We already know that $(\overline{\pi}^*K^*_Y\otimes K_{\ZZ})\vert_Y$ is 
non-canonically isomorphic to $H^{-2}$, therefore $\co_{\ZZ}(-2)$ is an 
extension of $(\overline{\pi}^*\ck^*_Y\otimes\ck_{\ZZ})\vert_Y$ over 
the whole $\ZZ$ that is,
\[\ext^{3-k}(\ZZ;j_*\co_Y,\overline{\pi}^*\ck^*_Y\otimes\ck_{\ZZ})\cong 
\ext^{3-k}(\ZZ ;j_*\co_Y,\co_{\ZZ}(-2)).\] 
Therefore, taking into account that $H^2(Y;\ck_Y)\cong 0$ and 
$H^1(Y;\ck_Y)\cong\C$ we conclude that
\[\ext^1(\ZZ ;j_*\co_Y,\co_{\ZZ}(-2))\cong 0\:\:\:\:\:{\rm 
and}\:\:\:\:\:\ext^2(\ZZ ;j_*\co_Y,\co_{\ZZ}(-2))\cong\C .\]
Referring to (\ref{ext-H-atteres}) we know 
$\ext^i(\ZZ ;\co_{\ZZ} ,\co_{\ZZ} (-2))\cong H^i(\ZZ
;\co_{\ZZ}(-2))$ consequently Lemma \ref{kohomologia} additionally 
yields
\[\ext^1(\ZZ ;\co_{\ZZ} ,\co_{\ZZ} (-2))\cong\C\:\:\:\:\:{\rm 
and}\:\:\:\:\:\ext^2(\ZZ ;\co_{\ZZ} ,\co_{\ZZ} (-2))\cong 0.\] 

Now consider the short exact sequence of coherent sheaves over $\ZZ$:
\[0\longrightarrow\cj_Y\longrightarrow\co_{\ZZ}\longrightarrow j_*\co_Y
\longrightarrow 0.\]
A segment of the associated long exact sequence of the global Ext groups 
with $\co_{\ZZ}(-2)$ looks like
\[\begin{array}{llll}
\hfill{\dots} 
& \longrightarrow \ext^1(\ZZ ;j_*\co_Y,\co_{\ZZ} (-2))
& \longrightarrow &\ext^1(\ZZ ;\co_{\ZZ} ,\co_{\ZZ} (-2))\longrightarrow\\
& \longrightarrow\ext^1(\ZZ ; \cj_Y,\co_{\ZZ} (-2))
& \stackrel{\delta}{\longrightarrow} & \ext^2(\ZZ ;j_*\co_Y,\co_{\ZZ} (-2))
\longrightarrow\\
& \longrightarrow\ext^2(\ZZ ;\co_{\ZZ} ,\co_{\ZZ} (-2)) 
& \longrightarrow 
& \dots\\
\end{array}\]
which gives
\begin{equation}
0\longrightarrow\C\longrightarrow \ext^1(\ZZ ; \cj_Y,\co_{\ZZ} 
(-2))\longrightarrow\C\longrightarrow 0
\label{hasitas}
\end{equation}
providing the result. $\Diamond$
\vspace{0.1in}

\noindent After these algebro-geometric preliminaries the time has come 
to relate our considerations so far with anti-instantons over a 
multi-Taub--NUT space. 

First we establish a link with harmonic functions. The embedding 
$i:Z\setminus Y\subset\ZZ$ 
together with (\ref{ext-H-atteres}) implies a homomorphism
\[i^*:\ext^1(\ZZ ;\cj_Y,\co_{\ZZ}(-2))\longrightarrow \ext^1(Z\setminus 
Y;\cj_Y,\co_{Z\setminus Y}(-2))\cong H^1(Z\setminus Y;\co_{Z\setminus 
Y}(-2)).\]
However for a fixed real line $Y\subset Z$ representing the point $y\in 
M_V$ via Penrose transform we have an isomorphism 
\[T:H^1(Z\setminus Y;\co_{Z\setminus Y}(-2))\cong{\rm
Ker}\:\triangle_{g_V}\vert_{M_V\setminus\{ y\}}.\]
Consequently one can think of the elements of 
$\ext^1(\ZZ ;\cj_Y,\co_{\ZZ}(-2))$ as 
harmonic functions. Therefore Lemma \ref{ext-lemma} provides us 
with a distinguished $2$-parameter family of complex-valued harmonic 
functions on $M_V\setminus\{ y\}$. These functions 
can easily be found by observing that the 
positive minimal Green function with singularity in $y\in M_V$ 
of Defintion \ref{green} provides us with a splitting of $\ext^1(\ZZ 
;\cj_Y,\co_{\ZZ}(-2))\cong\C^2$ as a vector 
space into two one dimensional summands. More precisely 
(\ref{hasitas}) gives 
\[0\longrightarrow\ext^1(\ZZ ;\co_{\ZZ} ,\co_{\ZZ} (-2))\longrightarrow
\ext^1(\ZZ ; \cj_Y,\co_{\ZZ} (-2))\stackrel{\delta}{\longrightarrow}
\ext^2(\ZZ ;j_*\co_Y,\co_{\ZZ} (-2))\longrightarrow 0\]
and by \cite[Theorem 1]{ati} the Green function defines a map 
\[\ext^2(\ZZ ;j_*\co_Y,\co_{\ZZ} (-2)) \longrightarrow \ext^1(\ZZ ; 
\cj_Y,\co_{\ZZ} (-2)),\]
so that all elements of the middle term canonically can be 
written in the form $\alpha +\beta$ with the first factor $\alpha$ 
coming from $\ext^1(\ZZ ;\co_{\ZZ} ,\co_{\ZZ} (-2))$ which is independent of 
$Y$ hence corresponds to a {\it constant} harmonic function $\mu\in\C$ 
on $(M_V, g_V)$ while $\beta =\lambda G(\cdot,y)$ with $\lambda\in\C$. 

Consequently we obtain the following (also cf. \cite{ati}):
\begin{theorem} 
Fix a point $y\in M_V$ of the multi-Taub--NUT space $(M_V, g_V)$ and 
consider the associated real line $Y\subset Z$ and the group $\ext^1(\ZZ 
;\cj_Y,\co_{\ZZ}(-2))\cong\C^2$. Then any element of this group can 
canonically be written in the form $\mu +\lambda G(\cdot ,y)$ where 
$G(\cdot ,y)$ is the unique minimal positive Green function concentrated at 
$y\in M_V$ and $\mu ,\lambda\in\C$ are constants. $\Diamond$
\label{harmonikus-ext}
\end{theorem}

\begin{remark}\rm Notice that in fact the above theorem holds for any
anti-half-conformally flat non-parabolic manifold whose twistor space 
admits a compactification satisfying Lemma \ref{kohomologia}.
\end{remark}

\noindent One can then see that the positive harmonic
functions of (\ref{dekompozicio})---from which our anti-instantons 
$[\nabla_{\tilde{A}_{y,\lambda}}]$ stem---can also be regarded as 
elements of $\ext^1(\ZZ ;\cj_Y,\co_{\ZZ}(-2))$ hence they describe 
extensions of two coherent sheaves
\[0\longrightarrow\co_{\ZZ}(-2)\longrightarrow\overline{\cf}(-1)
\longrightarrow\cj_Y\longrightarrow 0.\]
Restricting this to $Z$, the resulting sheaves $\cf (-1)$ will be 
locally free hence provide us with rank 2 holomorphic vector bundles 
$F(-1)$ over $Z$. 

Our aim is now to understand those twisted vector bundles together with 
their canonical sections among the aforementioned vector bundles over $Z$ which 
are associated to pure Green functions. We want to show that for two 
points $y,{\rm e}^{\ii\tau }y\in M_V$ on the same $S^1$-orbit, these 
twisted vector bundles are isomorphic. 
For this purpose we have to construct the real structure on $Z$ which 
gives rise to the multi-Taub--NUT geometry. It is sufficient to work 
over the open twistor space $Z\subset\ZZ$, because we are only interested in 
the real lines lying over points $y\in M_V$ at finite distance. Therefore we 
take a closer look at this space.

Restricting our construction in Lemma \ref{Zkompaktifikalas2} and in 
particular the singular fibration (\ref{szing.tvisztor.fibralas}) to the
finite part we obtain a model for the non-compact 
singular twistor space $\pi :Z^*\rightarrow\C P^1$ as the hypersurface 
$P_c=0$ where $P_c:(H^s\otimes L^c)\oplus (H^s\otimes L^{-c})\oplus 
H^2\rightarrow H^{2s}$ takes the shape
\begin{equation}
P_c(x,y,z):=(x\cdot 
\eta^c)(y\cdot\eta^{-c})-\prod\limits_{j=1}^s(z-p_j).
\label{veges.tvisztor.fibralas}
\end{equation}
(For clarity we remark that in (\ref{veges.tvisztor.fibralas}) and from 
now on without introducing extra notation all the bundles $H^k$ over 
$\C P^1$ will be pulled back to $H^2$.) 

But first we make a digression on the choice $k=l=s$ and $m=2$ in Lemmata 
\ref{Zkompaktifikalas1} and \ref{Zkompaktifikalas2} and note that 
it is also dictated by the requirement that the normal bundle of a real line 
$Y\subset Z$ must be isomorphic to $H\oplus H$ and $Z$ must possess a 
non-trivial real structure leading to the multi-Taub--NUT geometry. 

Indeed, to describe the normal bundle of a real 
line $Y\subset Z$ is equivalent to describe that of a generic real line 
$Y\subset Z^*$. Since it is a section of the bundle $\pi : 
Z^*\rightarrow\C P^1$ its normal bundle is the restriction 
of the tangent bundle of $H^s\oplus H^s\oplus H^2$ to a generic 
(non-singular) line i.e., the kernel of the map $(P_x, P_y, 
P_z):T(H^s\oplus H^s\oplus H^2)\rightarrow TH^{2s}$. This shows that 
$\langle c_1(N_Y), [Y]\rangle =s+s+2-2s=2$ moreover the type of the 
normal bundle is stable under small perturbations yielding $N_Y\cong 
H\oplus H$. An important consequence is that since $K_Z\vert_Y\cong 
K_Y\otimes\Lambda^2N^*_Y$ we therefore find $K_Z\cong\pi^*H^{-4}$ for 
the canonical bundle of $Z$ as we already noted in Sect. \ref{two}. 
Notice that this relation does not extend to the compactification $\ZZ$.
Indeed, for every real line $Y$ at finite distance the relation
$N_Y\cong H\oplus H$ continues to hold, but because the fibers of
the map (\ref{tvisztor.fibralas}) become compact, in $K_{\ZZ}$ there 
will be a twist coming from the canonical bundle of $\overline{X}$ as 
well. This explains the asymmetry between the middle cohomology groups 
in Lemma \ref{kohomologia}.

Next we move to the construction of the real 
structure on $Z$ following \cite{ete} --- originally due 
to \cite{bes,hit1,hit3} --- which gives rise to the multi-Taub--NUT geometry. 
As we have seen in Sect. \ref{two} the real structure on our twistor space 
must be induced by the antipodal map on $\C P^1$. Therefore, assuming 
$c\in\R$, we use (\ref{valos}) to construct induced maps 
\[r_c :H^0\left( H^2;\ce_{H^2}(H^s\otimes 
L^c)\right)\longrightarrow H^0\left( H^2;\ce_{H^2}(H^s\otimes L^{-c})\right)\] 
satisfying $r_c^2=(-1)^s{\rm Id}$. 
The {\it real structure $\tau_c:Z^*\rightarrow 
Z^*$ on the twistor space} given by $P_c=0$ in 
(\ref{veges.tvisztor.fibralas}) is defined to be
\begin{equation}
\tau_c(x,y,z):=\left( (-1)^s 
r_c(y\cdot\eta^{-c}),r_c(x\cdot\eta^c),-r(z)\right) .
\label{tauc}
\end{equation} 
A twistor line is called {\it real} if it is invariant under this map. 
We claim that the space $Z^*$ constructed in Lemma \ref{Zkompaktifikalas2}
together with this $\tau_c$ is indeed (a singular model of) the twistor 
space of the multi-Taub--NUT geometry. The simplest way of demonstrating this 
would be the derivation of the multi-Taub--NUT metric (\ref{metrika}) 
via Penrose' non-linear graviton construction outlined in Sect. 2. 
However, after a considerable hesitation we decided {\it not} to present 
this calculation here because it is quite long moreover it is already 
available in the literature for a long time (\cite[pp. 393-395]{bes}, 
\cite{hit1, hit3}).

We rather find the corresponding real lines by a factorization method 
(\cite[pp. 393-395]{bes} and \cite{hit1}). Using the notations and 
constructions of Lemma \ref{Zkompaktifikalas2} if we encode the 
NUT $q_j\in\R^3$ as the real section of 
$H^2$ whose shape over $U_b\subset\C P^1$ 
looks like $p_j([\frac{a}{b}:1])\vert_{U_b}
=\alpha_j\left(\frac{a}{b}\right)^2+2\beta_j\left(\frac{a}{b}\right)
-\overline{\alpha}_j$ as well as write a general point in $\R^3$ as $\zeta 
([\frac{a}{b}:1])\vert_{U_b} 
=\alpha\left(\frac{a}{b}\right)^2+2\beta\left(\frac{a}{b}\right)
-\overline{\alpha}$ then the roots 
$\rho_j,\sigma_j\in\C\subset\C P^1$ of the equation $\zeta -p_j=0$ are
\[\rho_j:=\frac{-(\beta -\beta_j)-\sqrt{(\beta 
-\beta_j)^2+\vert\alpha -\alpha_j\vert^2}}{\alpha 
-\alpha_j},\:\:\:\:\:\sigma_j:=\frac{-(\beta -\beta_j)+\sqrt{(\beta
-\beta_j)^2+\vert\alpha -\alpha_j\vert^2}}{\alpha -\alpha_j}.\] 
The real lines appear by simply factorizing in a 
$\tau_c$-invariant way the equation (\ref{veges.tvisztor.fibralas}) and 
imposing the reality condition (this already has been done
for $\zeta$). They are of the form
\[\left\{\begin{array}{ll}
\xi (\zeta ,[a:b])= & A\eta^c 
(\zeta ,[a:b])\prod\limits_{j=1}^s(a-\rho_jb)\\
\upsilon (\zeta ,[a:b])= & B\eta^{-c}(\zeta , [a:b])
\prod\limits_{j=1}^s(a-\sigma_jb)\\
\zeta ([a:b])= & \alpha a^2+2\beta ab -\overline{\alpha}b^2.
\end{array}\right.\]
To see that this is indeed a $\tau_c$-invariant factorization, 
the only non-trivial fact we have to check is that $\eta^{-c}$ goes to 
$\overline{\eta^c}$ if $c\in\R$ and {\it vice versa} on both charts 
$U'_a$ and $U'_b$. However for example for $[a:b]\in U'_a$ we have
$[-\overline{b}:\overline{a}]\in U'_b$, and $\eta^{-c}(\zeta
,[-\overline{b}:\overline{a}])={\rm
exp}(\overline{c(\overline{\alpha}(\frac{b}{a})-\beta )})$ which is
just $\overline{\eta^c(\zeta ,[a:b])}$ written on the chart $U'_a$.
The coefficients $A,B$ here are arbitrary constants 
satisfying $AB=\prod (\alpha - \alpha_i)$, the leading coefficient of
$\prod (\zeta -p_i)$. As this latter number is non-zero for 
$\alpha\not=\alpha_j$, say $B$ can be expressed in terms of $A$ and 
$\alpha$. But the reality constraint moreover implies
\[ \vert A\vert^2=\prod\limits_{j=1}^s\left( (\beta -\beta_j)+\sqrt{(\beta
-\beta_j)^2+\vert\alpha -\alpha_j\vert^2}\right) .\]
Hence the set of these real lines is locally parameterized\footnote{Note 
that the way of assigning the roots of $\zeta -p_j$ to $\xi$ and 
$\upsilon$ is well-defined only up to the action of the Galois group of 
the problem which explains that we have a local parameterization only.}
by ${\rm arg}A\in S^1$ and $({\rm Re}\:\alpha 
,{\rm Im}\:\alpha ,\beta)\in\R^3$ and provides us with a local chart of 
$M_V$ as it was constructed in Sect. \ref{three}. The choice 
$c=0$ gives rise to the ALE $A_{s-1}$ geometries while $c>0$ 
real provides us with the ALF $A_{s-1}$ i.e., the multi-Taub--NUT spaces. 
Hence from now on we will restrict to $c=1$ in accord with previous 
sections. A given real line will be denoted as $Y_{\alpha ,\beta ,A}$.

The next step is to express the transition matrices for the vector bundle 
$F_{\alpha ,\beta ,A}(-1)$ on $Z$ as well as its
section $s_{\alpha ,\beta ,A}\in H^0\left( Z;\co_Z
(F_{\alpha ,\beta ,A}(-1))\right)$, associated to the canonical element 
\[G(\cdot ,y_{\alpha ,\beta ,A})\in\ext^1(Z\setminus 
Y_{\alpha ,\beta ,A}\:;\:\cj_{Y_{\alpha ,\beta ,A}},\co_{Z\setminus 
Y_{\alpha ,\beta ,A}}(-2))\] 
representing the Green function. Here, we follow again \cite{ati} and 
first work over $Z^*$. Set
\[\begin{array}{ll} 
f := & x\cdot\eta^1 - \xi \\
g := & y\cdot\eta^{-1} - \upsilon \\ 
h := & z -\zeta . 
         \end{array}\] 
By definition, the real twistor line $Y_{\alpha ,\beta ,A}$ is given by 
the equations 
\[f = g = h = 0.\]
Define open subsets $U_{\xi}$ and $U_{\upsilon}$ by $\xi\not= 0$ and 
$\upsilon\not= 0$ respectively. Then, restricted to 
$U_{\xi}\cap Z^*$ the real line 
is given by the complete intersection $f=h=0$ as well as restricted to 
$U_{\upsilon}\cap Z^*$ it is the complete intersection $g=h=0$. 
Furthermore, for all generic $\zeta$ the polynomials $\zeta -p_i$ 
and $\zeta -p_j$ do not have common roots for $i\not= j$, so $U_{\xi}$ 
and $U_{\upsilon}$ cover $Z^*$. We restrict to such choices 
$\alpha ,\beta ,A$, but by continuity our results continue to hold in the 
general case too. Let us now set 
\[\theta := \frac{(x\cdot\eta^1)(y\cdot\eta^{-1})-\xi\upsilon}{h}=
\frac{xy-\xi\upsilon}{h},\]
a polynomial in $a,b$. Straightforward verification yields the identities 
\[f = -\frac{x\cdot\eta^1}{\upsilon} g + \frac{\theta}{\upsilon} h
\:\:\:\:\:{\rm and}\:\:\:\:\: 
    h = -\frac{h}{\upsilon} g + \frac{y\cdot\eta^{-1}}{\upsilon} h.\]
Obviously $f,g$ and $h$ are well-defined sections, vanishing on 
$Y_{\alpha ,\beta , A}$, of $H^s\otimes L^1$, $H^s\otimes L^{-1}$ and 
$H^2$ respectively. 

Hence, we can define a rank $2$ vector bundle 
$F_{\alpha ,\beta ,A}(-1)$ on $Z^*$ by gluing the sections $f$ and $h$ 
of the bundle $(H^s\otimes L^1) \oplus H^2$ restricted to $U_{\xi}$ 
with the sections $g$ and $h$ of the bundle $(H^s\otimes 
L^{-1}) \oplus H^2$ on $U_{\upsilon}$ using the gluing matrix 
\begin{equation}
   M_{\alpha ,\beta, A}:=\frac{1}{\upsilon}
  \begin{pmatrix}
       -x\cdot\eta^1 & \theta \\
        -h & y\cdot\eta^{-1}
  \end{pmatrix}
\label{gluingmatrix}
\end{equation}
on $U_{\xi}\cap U_{\upsilon}$. Moreover, we get a section 
$s_{\alpha ,\beta ,A}$ of this vector bundle by setting it to be equal to 
$(f,h)$ on $U_{\xi}$ and to $(g,h)$ on $U_{\upsilon}$. This 
section vanishes precisely on $Y_{\alpha ,\beta ,A}$. This bundle 
provides us with a Green function hence is positive \cite{gri,heb,var} 
moreover it is the minimal positive one since from \cite[Eq. 7.11]{ati} one 
can check that it converges to zero at infinity. 

Recall also that $S^1$ acts on $M_V$ as identified with 
the space of real sections by multiplication on the component $A$, or 
equivalently in view of the identity $AB=\prod (\alpha - \alpha_i)$, 
by inverse multiplication on $B$. Our aim is to show that if we consider the 
twistor line associated to a point 
$(\alpha ,\beta ,{\rm e}^{\ii\tau}A)$, then the vector bundle 
$F(-1)$ constructed for this point will be 
isomorphic to the one corresponding to $(\alpha ,\beta ,A)$. 
But the gluing matrix (\ref{gluingmatrix}) for these former data 
is related to that of the latter ones via
\[M_{\alpha ,\beta ,{\rm e}^{\ii\tau}A}={\rm 
e}^{-\ii\tau}M_{\alpha ,\beta ,A}.\] 
It follows that the change of trivializations $(f,h)\mapsto 
({\rm e}^{\ii\tau}f,{\rm e}^{\ii\tau}h)$ 
on $U_{\xi}$ induces an isomorphism between the bundles 
$F_{\alpha ,\beta ,A}(-1)$ and $F_{\alpha ,\beta ,{\rm e}^{\ii\tau}A}(-1)$. 

As we have seen we obtain $Z$ from the singular space $Z^*$ by blowing 
it up in finitely many points. Carrying out this procedure and then 
pulling back the bundles just constructed we come up with isomorphic 
vector bundles $F_{\alpha ,\beta ,A}(-1)$ and
$F_{\alpha ,\beta ,{\rm e}^{\ii\tau}A}(-1)$ over the regular 
(and non-compact) twistor space $Z$.

Consider the untwisted vector bundle $F_{\alpha ,\beta ,A}=F_{\alpha 
,\beta ,A}(-1)\otimes H$. First of all, it is clear that it is a 
holomorphic rank 2 vector bundle over $Z$. 

Secondly, the restriction of $F_{\alpha ,\beta ,A}$ to any real line in 
$Z$ is holomorphically trivial. This follows since the bundle we have 
constructed just corresponds to the Green function and this function is 
everywhere positive as we already noted above.

Thirdly, the anti-linear map (\ref{tauc}) induces a real structure on
$F_{\alpha ,\beta ,A}(-1)$ over $\tau_1$ simply because 
$\tau_1$ switches $H^s\otimes L^1$ and $H^s\otimes L^{-1}$ while maps 
$H^2$ to $H^2$. Therefore, since over $\C P^1$ the bundle $H$ carries a 
symplectic structure \cite{har}, the untwisted bundle $F_{\alpha ,\beta 
,A}= F_{\alpha ,\beta ,A}(-1)\otimes H$ also carries a {\it symplectic} 
structure 
$\tilde{\tau}_1: F_{\alpha ,\beta ,A}\rightarrow F_{\alpha ,\beta ,A}$ 
lying over $\tau_1$.

Since the same conclusions are also true for $F_{\alpha ,\beta
,{\rm e}^{\ii\tau}A}$ it then follows
that the corresponding untwisted vector bundles satisfy
the three conditions of the Atiyah--Ward correspondence (cf. the summary 
of the Atiyah--Ward correspondence in Sect. \ref{two}) hence they indeed 
give rise to anti-instantons. Since $F_{\alpha ,\beta ,A}$ and 
$F_{\alpha ,\beta ,{\rm e}^{\ii\tau}A}$ are isomorphic, we therefore 
conclude that the limiting anti-instantons are gauge equivalent i.e.,
$[\nabla_{\tilde{A}_{y,+\infty}}]=
[\nabla_{\tilde{A}_{{\rm e}^{\ii\tau}y,+\infty}}]$ as desired.


\section{Conclusion}


In this paper we considered SU(2) anti-instanton moduli spaces over the 
multi-Taub--NUT spaces containing solutions with integer energy. 
As a consequence of the general theory \cite{ete-jar} in principle 
anti-instantons with the same nice properties 
but with fractional energy also may exist over these spaces. The 
tantalazing question therefore arises whether or not these other moduli 
spaces are empty.

\end{document}